\documentclass[twoside,11pt]{article}
\usepackage{jnm}
\usepackage{paralist}
\usepackage{pifont}
\usepackage{amsmath,amsfonts,amssymb,cases}
\usepackage{empheq}
\usepackage{array,dcolumn,booktabs}
\usepackage[dvipsnames,svgnames]{xcolor}
\usepackage{tikz,todonotes}
\usepackage{subfig}
\usepackage{url}
\usepackage{microtype}
\usepackage{marvosym}
\usepackage{graphicx}
\usepackage{hyperref}
\usepackage{epstopdf}
\graphicspath{{images/}}
\usepackage{listings}
\newcommand{\ds}{\displaystyle}
\newcommand{\bg}{\begin{equation}}
\newcommand{\ed}{\end{equation}}
\usepackage{xspace}
\makeatletter
\DeclareRobustCommand{\freefem}
{\valign{\vfil\hbox{##}\vfil\cr
   \textsf{FreeFem\kern-.1em}\cr
   $\hbox{\fontsize{\sf@size}{0}\textbf{+\kern-0.05em+}}$\cr}\xspace%
}
\hypersetup{
backref=true, 
pagebackref=true,
hyperindex=true, 
colorlinks=true, 
breaklinks=true, 
urlcolor= black, 
linkcolor= black, 
citecolor = blue, 
 bookmarks=true, 
bookmarksopen=true, 
pdftitle={Solution of 2D Boussinesq systems with FreeFem++: the flat bottom case}, 
pdfauthor={Georges SADAKA}, 
pdfsubject={boussinesq_freefem_article} 
} 
\makeatother
\theoremstyle{plain}

\begin{document}
\date{}
\title{ {\sc \textbf{Solution of 2D Boussinesq systems with FreeFem++: the flat bottom case}}}
\author{{\sc \Large{Georges Sadaka}}\footnote{Universit\'e de Picardie Jules Verne, LAMFA CNRS UMR 7352, 33, rue Saint-Leu, 80039 Amiens, France, \url{http://lamfa.u-picardie.fr/sadaka/}  \Email\xspace \small{georges.sadaka@u-picardie.fr}}}
\maketitle
\textbf{Abstract} - We consider here different family of Boussinesq systems in two space dimensions. These systems approximate the three-dimensional Euler equations and consist of three coupled nonlinear dispersive wave equations that describe propagation of long surface waves of small amplitude in ideal fluids over a horizontal bottom and which was studied in \cite{Chen*09,Dou*07,Dou*09}\label{Dou*070}\label{Dou*090}\label{Chen*090}. We present here a \freefem code aimed at solving numerically these systems where a discretization using $\mathbb{P}_1$ finite element for these systems was taken in space and a second order Runge-Kutta scheme in time. We give the detail of our code where we use a mesh adaptation technique. An optimization of the used algorithm is done and a comparison of the solution for different Boussinesq family is done too. The results we obtained agree with those of the literature.

\begin{flushleft}

\end{flushleft}

\section{Introduction}
\begin{flushleft}
It has often been observed that variations of the bottom could influence the damping of the waves including extreme ones as {\it Tsunamis}: the coral reef or the underwater forests in the first shoreline, mangroves; these underwater reefs are also used to prevent corrosion effects of coastal (see P. Azerad {\it et al.} \cite{Aze*080} and \cite{Aze**080}\label{Aze*0801}\label{Aze**0801}). In these cases, the underwater relief damped the wave energy, in contrast, in other situations we seek to harness this energy: some companies even offer projects underwater reefs for erectile produce energy from waves (see \url{http://www.aquamarinepower.com/}).\vspace{.5cm}\\
In a continuum approximation, waves on the surface of an ideal fluid under the force of gravity are governed by the Euler equations. These are expected to provide a good model of irrotational waves on the surface of water, say, in situations where dissipative and surface tension effects may be safely ignored, see \cite{Mir*05}\label{Mir*050}.\vspace{.5cm}\\
In this paper, the attention is given to a multi-dimensional Boussinesq systems which describes approximately the propagation of small amplitude and long wave-length surface waves in a three-dimensional wave tank filled with an irrotational, incompressible and inviscid liquid under the influence of gravity, moving and/or variable bottom and surface pressure. The full Boussinesq systems in 2D has been derived in \cite{Sad11}\label{Sad111}.\vspace{.5cm}\\
Chen, Goubet, Dougalis, Mitsotakis and Saut have considered 2D models of Boussinesq systems with flat bottom \cite{Chen*09,Dou*07}\label{Dou*071}\label{Chen*091} then with variable bottom \cite{Chen09,Dut*07,Mit09}\label{Mit090}\label{Chen090}\label{Dut*070}, on the other hand Dutykh, Katsaounis and Mitsotakis have developed a code in finite volumes for the Boussinesq systems with variable bottom in 1D (\cite{Dut*110}\label{Dut*1101}) and Mitsotakis {\it et al.} in Galerkin finite elements (using B-splines \cite{Dou*07}\label{Dou*072}).\vspace{.5cm}\\
In order to solve numerically the Boussinesq system, we will use \freefem which is an open source platform to solve partial differential equations numerically, based on finite element methods. The \freefem platform has been developed to facilitate teaching and basic research through prototyping. For the moment this platform is restricted to the numerical simulations of problems which admit a variational formulation.\vspace{.5cm}\\
Thus, we develop a \freefem code for the simulation of Boussinesq equations with flat bottom, we first check that the simulations provided by our numerical code are consistent with the results of the recent literature, including the work of Dougalis, Mitsotakis et Saut \cite{Dou*07,Dou*09,Dou*10}\label{Dou*073}\label{Dou*091}\label{Dou*1004}. This establishes the adequacy of the chosen finite element discretization.\vspace{.5cm}\\
The article is organized as follows: first we discretize the problem in space by using finite element method and in time by using an explicit second order Runge-Kutta scheme, then we develop all the steps of the \freefem code to solve the problem by using the technique of mesh adaptation and at the end we present some numerical results where an optimization of the used algorithm is done and a comparison of the solution for different Boussinesq family is done too.
\end{flushleft}
\section{The Problem}
\subsection{Problem settings}\label{probset}
The 2D Boussineq system describe the surface wave propagation of small amplitude and large wave length. When considering an incompressible fluid flows in $\Omega\subset\mathbb{R}^2$, they are expressed as (see \cite{Mir*05}\label{Mir*051}):
\bg\label{BOUS}\begin{array}{rcl}\eta_t+\nabla\cdot  \textbf{V}+\nabla\cdot(\eta \textbf{V})+a\Delta\nabla\cdot \textbf{V}-b\Delta\eta_t&=&0;\\ \\
\textbf{V}_t+\nabla\eta+\frac{1}{2}\nabla|\textbf{V}|^2+c\Delta\nabla\eta-d\Delta \textbf{V}_t&=&0,\end{array}\ed
The variables in (\ref{BOUS}) are non-dimensional and unscaled: $\textbf{X}=(x,y)\in\Omega$ and $t>0$ are proportional to position along the channel and time, respectively, $\eta=\eta(\textbf{X},t)$ is proportional to the deviation of the free surface from its rest position,
$\textbf{V}=\textbf{V}(\textbf{X},t)=\left(\begin{array}{c} u(\textbf{X},t) \\  v(\textbf{X},t)\end{array}\right)=(u,v)^T=(u;v)$ is proportional to the horizontal velocity of the fluid at some height, $\nabla\textbf{*}=\left(\begin{array}{c}\partial_x\textbf{*} \\ \partial_y\textbf{*}\end{array}\right)$ is the gradient, $\nabla\cdot \left(\begin{array}{c}\star \\ \textbf{*}\end{array}\right)=\partial_x \star+\partial_y \textbf{*}$ is the divergence and $\Delta \textbf{*} =\partial_{xx}\textbf{*} +\partial_{yy}\textbf{*}$ is the laplacian. The coefficients $a,b,c$ and $d$ are given by the following formulas:
\bg\label{coefBOUS}a=\frac{1}{2}\left(\theta^2-\frac{1}{3}\right)\nu,b=\frac{1}{2}\left(\theta^2-\frac{1}{3}\right)(1-\nu),c=\frac{1}{2}\left(1-\theta^2\right)\mu,d=\frac{1}{2}\left(1-\theta^2\right)(1-\mu)\ed
where $\nu,\mu$ are real constants and $0\leqslant \theta \leqslant 1$. \newline
We note that the dispersive constants  $a,b,c$ and $d$ satisfy the physical constraints (see \cite{Bona*04}\label{Bona*041} for detail):
\bg\label{condabcd}a+b+c+d=\frac{1}{3} \mbox{ and } c+d\geq 0\ed
We now list some of the different family of Boussinesq systems 2D in Table \ref{TABBOUS2D} of the form (\ref{BOUS}):\newline

\begin{table}[!htb]
\begin{tabular}{|>{\centering}m{4cm}|>{\centering}m{2.5cm}|>{\centering}m{.5cm}|>{\centering}m{1.5cm}|>{\centering}m{2cm}|}
\hline
System & $\theta^2$ & $\nu$ & $\mu$ & References \tabularnewline
\hline
BBM-BBM &$2/3$ & $0$ & $0$ &  \cite{Chen09,Chen*09,Dou*07,Dou*09,Dou*10}.\label{Dou*101} \label{Dou*074} \label{Chen091} \label{Chen*092}\label{Dou*092} \tabularnewline
\hline
Bona-Smith & $2/3\leq \theta^2\leq 1$ & $0$ & $\ds\frac{4-6\theta^2}{3(1-\theta^2)}$ & \cite{Chen*09,Dou*07,Dou*09,Dou*10}. \label{Dou*102}\label{Chen*093}\label{Dou*075}\label{Dou*093}\tabularnewline
\hline
`` General '' Boussinesq & $0\leq \theta^2\leq 1$ & any & any & \cite{Chen*09,Dou*07,Dou*10}. \label{Dou*103}\label{Dou*076} \label{Chen*094}\tabularnewline
\hline
KdV-KdV & $2/3$ & $1$ & $1$ & \cite{Lin*113}\label{Lin*1131}\tabularnewline
\hline
\end{tabular}
\caption{\label{TABBOUS2D} Examples of Boussinesq systems in 2D.}
\end{table}
\begin{flushleft}
In \cite{Dou*07}\label{Dou*077}, V. Dougalis, D. Mitsotakis and J-C. Saut have studied the Well-Posedness of the Boussinesq systems (\ref{BOUS}) (where $b\neq 0$ and $d\neq 0$) and have shown that this system is at least nonlinearly well-posed locally; and in the case of KdV-KdV system (where $b=d=0,a=c=1/6$), F. Linares, D. Pilod and J-C. Saut proved recently in \cite{Lin*113}\label{Lin*1131} the Well-Posedness of this system.
\end{flushleft}
Following \cite{Wal*02}\label{Wal*021}, we can write (\ref{BOUS}) as:
\bg\label{BOUSN}\begin{array}{rcl}
\Upsilon-\Delta  \textbf{V} &=&0;\\ \\
\eta_t+\nabla\cdot  \textbf{V}+\nabla\cdot(\eta \textbf{V})+a\nabla\cdot\Upsilon-b\Delta\eta_t&=&0;\\ \\
\Theta-\Delta \eta &=&0;\\ \\
 \textbf{V}_t+\nabla\eta+\frac{1}{2}\nabla| \textbf{V}|^2+c\nabla\Theta-d\Delta \textbf{V}_t&=&0,\end{array}\ed
where $\Upsilon=\left(\Upsilon^1,\Upsilon^2\right)^T=\left(\Upsilon^1;\Upsilon^2\right)$.\newline

\section{Numerical Scheme}
In this section, we present the spatial discretization using finite element method with $\mathbb{P}_1$ continuous piecewise linear functions as shown in \cite{Wal*02}\label{Wal*022} and for the time marching scheme an explicit second order Runge-Kutta \cite{Dem91}\label{Dem911} scheme as used in \cite{Dou*07}\label{Dou*078}.\\
We will use in our code a mesh adaptation technic that we can use solving the problem by using the method based on the declaration of the problem obtained by the weak formulation of the system (\ref{BOUSN}); or by using the second method that consist to build matrices and vectors to solve the direct system $\mathbf{A} \mathbf{X}=\mathbf{B}$, where the matrix $\mathbf{A}$ and the vectors $\mathbf{X},\mathbf{B}$ will be defined in the sequel.
\subsection{Spatial discretization}
We let $\Omega$ be a convex, plane domain, let $\mathbf{T}_h$ denote a regular, quasi uniform triangulation of $\Omega$ with triangles of maximum size $h<1$ \cite{Bre*94}\label{Bre*941}, let $V_h=\{v_h\in C^0(\bar{\Omega}); v_h|_T \in \mathbb{P}_1(T),\forall T\in\mathbf{T}_h\}$ denote a finite-dimensional subspace of $H^1(\Omega)=\{u\in L^2(\Omega)\mbox{ s.t. } \frac{\partial u}{\partial x},\frac{\partial u}{\partial y} \in L^2(\Omega)\}$ where $\mathbb{P}_1$ is the set of polynomials of $ \mathbb{R}$ of degrees $\leq 1$ and let $\left\langle\cdot ; \cdot\right\rangle$ denote the $L^2$ inner product on $\Omega$.\\
Consider the weak formulation of the system (\ref{BOUSN}), find $\eta_h, u_h, v_h\in V_h$ such that $\forall \phi_h\in V_h$ we have:
\bg\label{BOUSWeak}\begin{array}{rcl} \left\langle\Upsilon^1_h;\phi_h\right\rangle-\left\langle \Delta u_h;\phi_h\right\rangle=0;\quad \left\langle\Upsilon^2_h;\phi_h\right\rangle - \left\langle \Delta v_h;\phi_h\right\rangle=0;\quad \left\langle\Theta_h;\phi_h\right\rangle-\left\langle \Delta \eta_h;\phi_h\right\rangle&=&0;
\\ \Big\langle(Id-b\Delta)\eta_{ht}+\nabla\cdot\left(u_h;v_h\right)+\eta_{hx} u_h+\eta_h u_{hx}+\eta_{hy} u_h+\eta_h u_{hy}+a\nabla\cdot\left(\Upsilon^1_h;\Upsilon^2_h\right);\phi_h\Big\rangle&= &0 ;
\\  \Big\langle (Id-d\Delta)u_{ht}+\eta_{hx}+ u_h u_{hx}+ v_h v_{hx}+c\Theta_{hx};\phi_h\Big\rangle&=&0;
\\ \Big\langle (Id-d\Delta)v_{ht}+\eta_{hy}+ u_h u_{hy}+ v_h v_{hy}+c\Theta_{hy};\phi_h\Big\rangle&=&0.\end{array}\ed
To simplify, we denote $\Phi=\phi_h,\mathbf{E}=\eta_h,\mathbf{U}= u_h,\mathbf{V}= v_h,\mathbf{T}=\Theta_h,\mathbf{P}=\Upsilon_h^1$ and $\mathbf{Q}=\Upsilon_h^2$, so the system (\ref{BOUSWeak}) is equivalent to the following system:
\bg\label{NBOUSWeak}\begin{array}{rcl} 
\multicolumn{3}{c}
{\left\langle\mathbf{P};\Phi\right\rangle\quad=\quad\left\langle\Delta\mathbf{U};\Phi\right\rangle;\quad\left\langle\mathbf{Q};\Phi\right\rangle\quad=\quad\left\langle\Delta\mathbf{V};\Phi\right\rangle;\quad \left\langle\mathbf{T};\Phi\right\rangle\quad=\quad\left\langle\Delta\mathbf{E};\Phi\right\rangle;}\\ \Big\langle(Id-b\Delta)\partial_t\mathbf{E};\Phi\Big\rangle&=&-\left\langle\nabla\cdot(\mathbf{U};\mathbf{V})+\mathbf{E}_x\mathbf{U}+\mathbf{E}\mathbf{U}_x+\mathbf{E}_y\mathbf{V}+\mathbf{E}\mathbf{V}_y+a\nabla\cdot\left(\mathbf{P};\mathbf{Q}\right);\Phi\right\rangle\\&=&-\mathbf{F}\left(\mathbf{E},\mathbf{U},\mathbf{V},\mathbf{P},\mathbf{Q}\right);\\ \Big\langle(Id-d\Delta)\partial_t\mathbf{U};\Phi\Big\rangle&=&-\left\langle \mathbf{E}_{x}+\mathbf{U}\mathbf{U}_{x}+\mathbf{V}\mathbf{V}_{x}+c\mathbf{T}_x;\Phi\right\rangle=-\mathbf{G}\left(\mathbf{E},\mathbf{U},\mathbf{V},\mathbf{T}\right);\\\Big\langle(Id-d\Delta)\partial_t\mathbf{V};\Phi\Big\rangle&=&-\left\langle \mathbf{E}_{y}+\mathbf{U}\mathbf{U}_{y}+\mathbf{V}\mathbf{V}_{y}+c\mathbf{T}_y;\Phi\right\rangle = -\mathbf{H}\left(\mathbf{E},\mathbf{U},\mathbf{V},\mathbf{T}\right).\end{array}\ed

\subsection{Time marching scheme}
Our method is based on an explicit second order Runge-Kutta scheme. To this end, let us denote by $(\mathbf{E}^{n+1},\mathbf{U}^{n+1},\mathbf{V}^{n+1})$ and $(\mathbf{E}^n,\mathbf{U}^n,\mathbf{V}^n,\mathbf{P}^n,\mathbf{Q}^n,\mathbf{T}^n)$ the approximate value at time $t=t^{n+1}$ and $t=t^n$, respectively and by $\delta t$ the time step size. Then, by using (\ref{NBOUSWeak}), the unknown fields at time $t=t^{n+1}$ are defined as the solution of the system
\bg\label{TBOUSW}\left\{\begin{array}{l}
\left\langle\mathbf{P}^n;\Phi\right\rangle=\left\langle\Delta\mathbf{U}^n;\Phi\right\rangle;\qquad
\left\langle\mathbf{Q}^n;\Phi\right\rangle=\left\langle\Delta\mathbf{V}^n;\Phi\right\rangle;\qquad
\left\langle\mathbf{T}^n;\Phi\right\rangle=\left\langle\Delta\mathbf{E}^n ;\Phi\right\rangle;\\
\langle\mathbf{E}^{n+1};\Phi\rangle=\langle \mathbf{E}^n+\ds\frac{\mathbf{E}^{k1}+\mathbf{E}^{k2}}{2};\Phi\rangle; \langle\mathbf{U}^{n+1};\Phi\rangle=\langle \mathbf{U}^n+\ds\frac{\mathbf{U}^{k1}+\mathbf{U}^{k2}}{2};\Phi\rangle; \\\langle\mathbf{V}^{n+1};\Phi\rangle=\langle \mathbf{V}^n+\ds\frac{\mathbf{V}^{k1}+\mathbf{V}^{k2}}{2};\Phi\rangle.\end{array}\right.\ed
where:
\bg\label{PSIUVK1}\left\{\begin{array}{l}\left\langle(Id-b\Delta)\mathbf{E}^{k1};\Phi\right\rangle=-\delta t\cdot\mathbf{F}\left(\mathbf{E}^n,\mathbf{U}^n,\mathbf{V}^n,\mathbf{P}^n,\mathbf{Q}^n\right);\\ \left\langle(Id-d\Delta)\mathbf{U}^{k1};\Phi\right\rangle=-\delta t\cdot\mathbf{G}\left(\mathbf{E}^n,\mathbf{U}^n,\mathbf{V}^n,\mathbf{T}^n\right);\\\left\langle(Id-d\Delta)\mathbf{V}^{k1};\Phi\right\rangle=-\delta t\cdot\mathbf{H}\left(\mathbf{E}^n,\mathbf{U}^n,\mathbf{V}^n,\mathbf{T}^n\right);\\
\left\langle\mathbf{P}^{k1};\Phi\right\rangle=\left\langle\mathbf{U}^{k1}_{xx}+\mathbf{U}^{k1}_{yy};\Phi\right\rangle;  \left\langle\mathbf{Q}^{k1};\Phi\right\rangle=\left\langle\mathbf{V}^{k1}_{xx}+\mathbf{V}^{k1}_{yy};\Phi\right\rangle;  \left\langle\mathbf{T}^{k1};\Phi\right\rangle=\left\langle\mathbf{E}_{xx}^{k1}+\mathbf{E}_{yy}^{k1} ;\Phi\right\rangle.\end{array}\right.\ed
and
\bg\label{PSIUVK2}\left\{\begin{array}{rcl} \left\langle(Id-b\Delta)\mathbf{E}^{k2};\Phi\right\rangle&=&-\delta t\cdot\mathbf{F}\left(\mathbf{E}^n+\mathbf{E}^{k1},\mathbf{U}^n+\mathbf{U}^{k1},\mathbf{V}^n+\mathbf{V}^{k1},\mathbf{P}^n+\mathbf{P}^{k1},\mathbf{Q}^n+\mathbf{Q}^{k1}\right);\\ \left\langle(Id-d\Delta)\mathbf{U}^{k2};\Phi\right\rangle&=&-\delta t\cdot\mathbf{G}\left(\mathbf{E}^n+\mathbf{E}^{k1},\mathbf{U}^n+\mathbf{U}^{k1},\mathbf{V}^n+\mathbf{V}^{k1},\mathbf{T}^n+\mathbf{T}^{k1}\right);\\\left\langle(Id-d\Delta)\mathbf{V}^{k2};\Phi\right\rangle&=&-\delta t\cdot\mathbf{H}\left(\mathbf{E}^n+\mathbf{E}^{k1},\mathbf{U}^n+\mathbf{U}^{k1},\mathbf{V}^n+\mathbf{V}^{k1},\mathbf{T}^n+\mathbf{T}^{k1}\right) .\end{array}\right.\ed
By integrating by parts where we have second order derivative and by developing all the terms of the first order derivative in (\ref{TBOUSW}), (\ref{PSIUVK1}) and (\ref{PSIUVK2}), we deduce:
\bg\label{BOUSWFINAL1}\begin{array}{c}
\left\langle\mathbf{P}^n;\Phi\right\rangle=-\left\langle\nabla\mathbf{U}^{n};\nabla\Phi\right\rangle+\left\langle\ds\frac{\partial\mathbf{U}^n}{\partial n};\Phi\right\rangle_{\partial \Omega};
\left\langle\mathbf{Q}^n;\Phi\right\rangle=-\left\langle\nabla\mathbf{V}^{n};\nabla\Phi\right\rangle+\left\langle\ds\frac{\partial\mathbf{V}^n}{\partial n};\Phi\right\rangle_{\partial \Omega};\\ \\
\left\langle\mathbf{T}^n;\Phi\right\rangle=-\left\langle\nabla\mathbf{E}^{n};\nabla\Phi\right\rangle+\left\langle\ds\frac{\partial\mathbf{E}^n}{\partial n};\Phi\right\rangle_{\partial \Omega};
\end{array}\ed

\bg\label{BOUSWFINAL2}\left\{\begin{array}{rcl}\left\langle\mathbf{E}^{k1};\Phi\right\rangle+b\left\langle\nabla\mathbf{E}^{k1};\nabla\Phi\right\rangle-b\left\langle\ds\frac{\partial\mathbf{E}^{k1}}{\partial n};\Phi\right\rangle_{\partial \Omega}&=&-\delta t\cdot\mathbf{F}\left(\mathbf{E}^n,\mathbf{U}^n,\mathbf{V}^n,\mathbf{P}^n,\mathbf{Q}^n\right);\\\left\langle\mathbf{U}^{k1};\Phi\right\rangle+d\left\langle\nabla\mathbf{U}^{k1};\nabla\Phi\right\rangle-d\left\langle\ds\frac{\partial\mathbf{U}^{k1}}{\partial n};\Phi\right\rangle_{\partial \Omega}&=&-\delta t\cdot\mathbf{G}\left(\mathbf{E}^n,\mathbf{U}^n,\mathbf{V}^n,\mathbf{T}^n\right);\\\left\langle\mathbf{V}^{k1};\Phi\right\rangle+d\left\langle\nabla\mathbf{V}^{k1};\nabla\Phi\right\rangle-d\left\langle\ds\frac{\partial\mathbf{V}^{k1}}{\partial n};\Phi\right\rangle_{\partial \Omega}&=&-\delta t\cdot\mathbf{H}\left(\mathbf{E}^n,\mathbf{U}^n,\mathbf{V}^n,\mathbf{T}^n\right);\\
\left\langle\mathbf{P} ^{k1};\Phi\right\rangle=-\left\langle\nabla\mathbf{U}^{k1};\nabla\Phi\right\rangle+\left\langle\ds\frac{\partial\mathbf{U} ^{k1}}{\partial n};\Phi\right\rangle_{\partial \Omega};& &\\
\left\langle\mathbf{Q} ^{k1};\Phi\right\rangle=-\left\langle\nabla\mathbf{V}^{k1};\nabla\Phi\right\rangle+\left\langle\ds\frac{\partial\mathbf{V} ^{k1}}{\partial n};\Phi\right\rangle_{\partial \Omega};& & \\
\left\langle\mathbf{T} ^{k1};\Phi\right\rangle=-\left\langle\nabla\mathbf{E}^{k1};\nabla\Phi\right\rangle+\left\langle\ds\frac{\partial\mathbf{E} ^{k1}}{\partial n};\Phi\right\rangle_{\partial \Omega};& &
\end{array}\right.\ed
and
\bg\label{BOUSWFINAL3}\left\{\begin{array}{l}\left\langle\mathbf{E}^{k2};\Phi\right\rangle+b\left\langle\nabla\mathbf{E}^{k2};\nabla\Phi\right\rangle-b\left\langle\ds\frac{\partial\mathbf{E}^{k2}}{\partial n};\Phi\right\rangle_{\partial \Omega}=\\
\qquad\qquad\qquad-\delta t\cdot\mathbf{F}\left(\mathbf{E}^n+\mathbf{E}^{k1},\mathbf{U}^n+\mathbf{U}^{k1},\mathbf{V}^n+\mathbf{V}^{k1},\mathbf{P}^n+\mathbf{P}^{k1},\mathbf{Q}^n+\mathbf{Q}^{k1}\right);\\ \left\langle\mathbf{U}^{k2};\Phi\right\rangle+d\left\langle\nabla\mathbf{U}^{k2};\nabla\Phi\right\rangle-d\left\langle\ds\frac{\partial\mathbf{U}^{k2}}{\partial n};\Phi\right\rangle_{\partial \Omega}=\\
\qquad\qquad\qquad-\delta t\cdot\mathbf{G}\left(\mathbf{E}^n+\mathbf{E}^{k1},\mathbf{U}^n+\mathbf{U}^{k1},\mathbf{V}^n+\mathbf{V}^{k1},\mathbf{T}^n+\mathbf{T}^{k1}\right);\\\left\langle\mathbf{V}^{k2};\Phi\right\rangle+d\left\langle\nabla\mathbf{V}^{k2};\nabla\Phi\right\rangle-d\left\langle\ds\frac{\partial\mathbf{V}^{k2}}{\partial n};\Phi\right\rangle_{\partial \Omega}=\\
\qquad\qquad\qquad-\delta t\cdot\mathbf{H}\left(\mathbf{E}^n+\mathbf{E}^{k1},\mathbf{U}^n+\mathbf{U}^{k1},\mathbf{V}^n+\mathbf{V}^{k1},\mathbf{T}^n+\mathbf{T}^{k1}\right) .\end{array}\right.\ed
{\bf Remark:} 
It's easy with \freefem to define boundary condition, in fact if we have the Dirichlet Boundary Conditions on a border $\Gamma_1\subset \mathbb{R}$ like $\mathbf{U}|_{\Gamma_1}=f$, then it is defined as {\ttfamily\textcolor{red}{on}(gamma1,u=f)}, where {\ttfamily{u}} is the unknown function in the problem. We note that the Neumann Boundary Conditions on $\Gamma_2\subset \mathbb{R}$, like $\ds\frac{\partial \mathbf{U}}{\partial n}|_{\Gamma_2}=g$, appear in the Weak formulation of the problem after integrating by parts for example in the system (\ref{BOUSWFINAL1}) we have $\left\langle\ds\frac{\partial\mathbf{U}}{\partial n};\Phi\right\rangle_{\Gamma_2}=\left\langle g;\Phi\right\rangle_{\Gamma_2}=\displaystyle\int_{\Gamma_2} g\cdot \Phi$  which is defined in \freefem by {\ttfamily\textcolor{red}{int1d}(Th,gamma2)(g*phi)} where  {\ttfamily Th} is the triangulated domain of $\Omega$. We will see in the next section how it's also easy to define the Bi-Periodic Boundary Conditions.\\
We remark also that the system (\ref{TBOUSW}) can be written on the following matrix form:
\bg\label{BOUSAXB}\underbrace{\left(\begin{array}{ccc}  \mathbf{M} &  0 & 0 \\ 0 & \mathbf{M} & 0 \\ 0 & 0 & \mathbf{M} \end{array}\right)}_{\mathbf{A}}\cdot\underbrace{\left(\begin{array}{c} \mathbf{E}^{n+1} \\ \mathbf{U}^{n+1} \\ \mathbf{V}^{n+1} \end{array}\right)}_{\mathbf{X}}=\underbrace{\left(\begin{array}{c} \langle \mathbf{E}^n+\ds\frac{\mathbf{E}^{k1}+\mathbf{E}^{k2}}{2};\Phi\rangle \\  \langle \mathbf{U}^n+\ds\frac{\mathbf{U}^{k1}+\mathbf{U}^{k2}}{2};\Phi\rangle \\  \langle \mathbf{V}^n+\ds\frac{\mathbf{V}^{k1}+\mathbf{V}^{k2}}{2};\Phi\rangle \end{array}\right)}_{\mathbf{B}} \ed
where $\mathbf{M}_{ij}=\ds\int_{\Omega}\phi_i\phi_j\mbox{d}x\mbox{d}y$ is the mass matrix.
\subsubsection*{Algorithm 1:}\label{ALGO1} $\,\,$Finally, to solve the systems (\ref{TBOUSW}), (\ref{BOUSWFINAL1}), (\ref{BOUSWFINAL2}) and (\ref{BOUSWFINAL3}), we follow as:
$$
\begin{array}{ll}
\mbox{\tt Set } & \mathbf{E}^n=\mathbf{E}^{0}=\eta_{h0}=\eta_{0}, \mathbf{U}^n=\mathbf{U}^0=u_{h0}=u_0, \mathbf{V}^n=\mathbf{V}^0=v_{h0}=v_0\\
\mbox{\tt Set } & \mathbf{P}^n=P_{h0}, \mathbf{Q}^n=Q_{h0}, \mathbf{T}^n=T_{h0},\mathbf{P}^{k1}=P_{hk1}, \mathbf{Q}^{k1}=Q_{hk1}, \mathbf{T}^{k1}=T_{hk1}\\
\mbox{\tt Set} & \mathbf{E}^{n+1}=\eta_{h}, \mathbf{U}^{n+1}=u_h, \mathbf{V}^{n+1}=v_h\\
\mbox{\tt Set} & \mathbf{E}^{k1}=\eta_{hk1}, \mathbf{U}^{k1}=u_{hk1}, \mathbf{V}^{k1}=v_{hk1},\mathbf{E}^{k2}=\eta_{hk2}, \mathbf{U}^{k2}=u_{hk2}, \mathbf{V}^{k2}=v_{hk2}\\
\mbox{\tt For $t=0:\delta t:T$ } & \\
&\mbox{\tt Mesh adaptation, (optional)} \\
&\mbox{\tt Compute } P_{h0},Q_{h0},T_{h0}\qquad \mbox{\tt Compute } \eta_{hk1}, u_{hk1}, v_{hk1}\\
&\mbox{\tt Compute } P_{hk1},Q_{hk1},T_{hk1}\qquad \mbox{\tt Compute } \eta_{hk2}, u_{hk2}, v_{hk2}\\
&\mbox{\tt Compute } \eta_h, u_h, v_h\\
\qquad\qquad\mbox{\tt Set} & \eta_{h0}=\eta_h, u_{h0}=u_h,v_{h0}=v_h\\
 \mbox{\tt End for}&
\end{array}
$$
\subsubsection*{Algorithm 2:}\label{ALGO2}
$\,\,$Another method to solve the systems (\ref{BOUSWFINAL1}), (\ref{BOUSWFINAL2}) and (\ref{BOUSWFINAL3}), taking into account (\ref{BOUSAXB}):
$$
\begin{array}{ll}
\mbox{\tt Set } & \mathbf{E}^n=\mathbf{E}^{0}=\eta_{h0}=\eta_{0}, \mathbf{U}^n=\mathbf{U}^0=u_{h0}=u_0, \mathbf{V}^n=\mathbf{V}^0=v_{h0}=v_0\\
\mbox{\tt Set } & \mathbf{P}^n=P_{h0}, \mathbf{Q}^n=Q_{h0}, \mathbf{T}^n=T_{h0},\mathbf{P}^{k1}=P_{hk1}, \mathbf{Q}^{k1}=Q_{hk1}, \mathbf{T}^{k1}=T_{hk1}\\
\mbox{\tt Set} & \mathbf{E}^{n+1}=\eta_{h}, \mathbf{U}^{n+1}=u_h, \mathbf{V}^{n+1}=v_h\\
\mbox{\tt Set} & \mathbf{E}^{k1}=\eta_{hk1}, \mathbf{U}^{k1}=u_{hk1}, \mathbf{V}^{k1}=v_{hk1},\mathbf{E}^{k2}=\eta_{hk2}, \mathbf{U}^{k2}=u_{hk2}, \mathbf{V}^{k2}=v_{hk2}\\
\mbox{\tt Compute } & \mathbf{A} \mbox{\tt (if we want to use the mesh adaptation }\\
& \mbox{\tt we must compute $\mathbf{A}$ in the for-loop time)}\\ 
\mbox{\tt For $t=0:\delta t:T$ } & \\
&\mbox{\tt Mesh adaptation, (optional)} \\
&\mbox{\tt Update } \eta_{h0}=\eta_{h0};u_{h0}=u_{h0};v_{h0}=v_{h0};\eta_h=\eta_h;u_h=u_h;v_h=v_h;\\
& \mbox{\tt (with mesh adaptation)}\\
&\mbox{\tt Compute } P_{h0},Q_{h0},T_{h0}\qquad \mbox{\tt Compute } \eta_{hk1}, u_{hk1}, v_{hk1}\\
&\mbox{\tt Compute } P_{hk1},Q_{hk1},T_{hk1}\qquad \mbox{\tt Compute } \eta_{hk2}, u_{hk2}, v_{hk2}\\
&\mbox{\tt Compute } \mathbf{A}\mbox{\tt (with mesh adaptation)}\\
&\mbox{\tt Set } \mathbf{X}=[\eta_h, u_h, v_h]\qquad \mbox{\tt Compute }  \mathbf{B}\qquad \mbox{\tt Solve }  \mathbf{A}\mathbf{X}=\mathbf{B}\\
\qquad\qquad\mbox{\tt Set} & \eta_{h0}=\eta_h, u_{h0}=u_h,v_{h0}=v_h\\
 \mbox{\tt End for}&
\end{array}
$$
\section{Code}
In this section we will present by details all the step of the \freefem code to solve (\ref{TBOUSW}) to (\ref{BOUSAXB}).
\subsection{Declaration of the problems}
Note that in \freefem the scalar product in $L^2$: $\left\langle  \textbf{\Large{.}},\phi_h\right\rangle=\displaystyle\int_\Omega \textbf{\Large{.}}\cdot \phi_h=$ {\ttfamily{\textcolor{red}{int2d}(Th)( \textbf{\Large{.}}*phih )} }; also we can define a macro for the right hand side function $\mathbf{F}\left(\mathbf{E},\mathbf{U},\mathbf{V},\mathbf{P},\mathbf{Q}\right)$, $\mathbf{G}\left(\mathbf{E},\mathbf{U},\mathbf{V},\mathbf{T}\right),\mathbf{H}\left(\mathbf{E},\mathbf{U},\mathbf{V},\mathbf{T}\right)$ defined in (\ref{NBOUSWeak}) using the keyword  \textcolor{red}{\ttfamily{macro}}, that will be used in the sequence, as:
\begin{lstlisting}[firstnumber=last]
macro F(e,u,v,p,q)(div(u,v)+dx(e)*u+e*dx(u)+dy(e)*v+e*dy(v)+a*div(p,q))//
macro G(e,u,v,t)(dx(e)+dx(u)*(u)+dx(v)*(v)+c*dx(t))//
macro H(e,u,v,t)(dy(e)+dy(u)*(u)+dy(v)*(v)+c*dy(t))//
\end{lstlisting}
We note that all the variable ({\ttfamily{e,u,v,p,q,t}}) used in the macro are dummies.\\
We declare the problem for $\mathbf{U}^{k1}$ defined is the system (\ref{BOUSWFINAL2}) and for $\mathbf{E}^{n+1}$ defined in the system (\ref{TBOUSW}) as:
\begin{lstlisting}[firstnumber=last]
problem UHK1(uhk1,phih) = int2d(Th)(uhk1*phih) + int2d(Th)(grad(uhk1)'*grad(phih)*d) + int2d(Th)( G(etah0,uh0,vh0,Th0)*phih*dt)+"Boundary Conditions of uh for uhk1";
problem ETAH(etah,phih) = int2d(Th)(etah*phih) - int2d(Th)(etah0*phih) - int2d(Th)((etahk1 + etahk2) * phih /2.);
\end{lstlisting}
The declaration of the problems for all other variables are written in the same form.\\

\textbf{Remark:} In order to make our code faster, we can use the keyword \textcolor{red}{\ttfamily{init}} in the declaration of the problem. When \textcolor{red}{\ttfamily{init}}{\ttfamily{=0}} the mass matrix is computed and when \textcolor{red}{\ttfamily{init}}{\ttfamily{=1}} the mass matrix is reused so it is much faster after the first iteration.

\subsection{Solve of the problems}
To solve all the problems defined above, we make a for-loop time and we call the problems by their names when we want them to be solved, then we update the data and at the end we plot the solution using the keyword \textcolor{red}{\ttfamily{plot}}.\newline
We note that in each iteration of the for-loop a mesh adaptation will be done which depend on the error ({\ttfamily{err}}) which is the $\mathbb{P}_1$ interpolation error level, 
where {\ttfamily{hmin}} is the minimum edge size and \textcolor{red}{\ttfamily{nbvx}} is the maximum number of vertices generated by the mesh generator.
\begin{lstlisting}[firstnumber=last]
for (real t=0.;t<=T;t+=dt){
	Th=adaptmesh(Th,etah0,err=1e-4,hmin=Dx,nbvx=1e6); // we can use adaptmesh each 10 iterations or more.
	PH0;		QH0;		TH0;
	ETAHK1;		UHK1;		VHK1;
	PHK1;		QHK1;		THK1;
	etah0pk1=etah0+etahk1;	uh0pk1=uh0+uhk1;	vh0pk1=vh0+vhk1;
	Th0pk1=Th0+Thk1;		Ph0pk1=Ph0+Phk1;	Qh0pk1=Qh0+Qhk1;
	ETAHK2;		UHK2;		VHK2;
	ETAH;		UH;			VH;
	etah0=etah;	uh0=uh;	vh0=vh; 	//update of the data
	plot(etah0,cmm="t="+t+"sec",fill=true,value=true,dim=3); 
}
\end{lstlisting}
In order to use the second method, we build the matrix $\mathbf{A}$ before the for-loop time as:
\begin{lstlisting}[firstnumber=last]
varf Mass(u, phih) = int2d(Th)( u * phih );		
matrix A, MASS;
MASS = Mass(Vh,Vh);
A = [[MASS, 0, 0],[0, MASS, 0],[0, 0, MASS]];
set(A,solver=GMRES); // to be set
\end{lstlisting}
Then we build the vector $\mathbf{B}$ in the for-loop time as:
\begin{lstlisting}[firstnumber=last]
for (real t=0.;t<=T;t+=dt){
PH0;		QH0;		TH0;
ETAHK1;		UHK1;		VHK1;
PHK1;		QHK1;		THK1;
etah0pk1=etah0+etahk1;	uh0pk1=uh0+uhk1;	vh0pk1=vh0+vhk1;
Th0pk1=Th0+Thk1;		Ph0pk1=Ph0+Phk1;	Qh0pk1=Qh0+Qhk1;	
ETAHK2;		UHK2;		VHK2;
Vh B1, B2, B3, etahk1pk2D2, uhk1pk2D2, vhk1pk2D2;
real[int] B(3*Vh.ndof), X(3*Vh.ndof), X0(3*Vh.ndof), W(3*Vh.ndof);
etahk1pk2D2 = .5*etahk1 + .5*etahk2;
uhk1pk2D2 = .5*uhk1 + .5*uhk2;
vhk1pk2D2 = .5*vhk1 + .5*vhk2;
X0=[etah0[], uh0[], vh0[]];
B1[]=MASS*etahk1pk2D2[];
B2[]=MASS*uhk1pk2D2[];
B3[]=MASS*vhk1pk2D2[];
B=[B1[],B2[],B3[]];
X = A^-1*B;
W = X + X0;
[etah[], uh[], vh[]] = W;
etah0=etah;	uh0=uh;	vh0=vh; 	//update of the data
plot(etah0,cmm="t="+t+"sec",fill=true,value=true,dim=3); 
}
\end{lstlisting}
Finally, if we want to use mesh adaptation in the second method, we must compute the matrix A in the for-loop time.

\subsection{Numerical simulations}
In the sequel, we present the results of numerical simulations of the evolution of initially localized heaps of fluid of initial velocity zero. Unless specified, all computations were performed on the square $\Omega= ]-40,40[\times ]-40,40[$, a $\mathbb{P}_1$ continuous piecewise linear functions was used for the finite element space and  for all the numerical simulations, we work with the space discretization $\Delta x=0.5$ and the time step $\Delta t=0.1$.\newline
\subsubsection{Rate of convergence}
At the beginning, we prove in the figure below, that the RK2 time scheme considered for the BBM-BBM system is of order 2. In this example, we took zero Dirichlet homogenous Boundary Conditions for $\eta_{h}$, $ u_{h}$ and $ v_{h}$ on the whole boundary and we have consider the following exacts solutions:
$$ \eta_{ex}=e^t\cdot \sin(\pi x)\cdot (y-1)\cdot y,$$ 
$$u_{ex}=e^t\cdot x \cdot \cos(3\pi x/2)\cdot \sin(\pi y),$$ 
$$v_{ex}=e^t\cdot\sin(\pi x)\cdot \cos(3\pi y/2)\cdot y.$$
Then, we compute the corresponding right hand side in order to obtain the $L^2$ norm of the error between the exact solution and the numerical one in the table below.
\begin{table}[!h]
\begin{center}
\begin{tabular}{|c|c|c|c|}
\hline
N & $|\eta_h-\eta_{ex}|_{L^2}$ & $|u_h-u_{ex}|_{L^2}$ & $|v_h-v_{ex}|_{L^2}$\\
\hline
10& 0.00871494 & 0.0233966& 0.0230945\\
\hline
20&  0.00265707&0.00641675&0.00632314\\
\hline
40& 0.000670301&0.00160223&0.00157848\\
\hline
80&0.0001817&0.000419198&0.000412791\\
\hline
160&4.80657e-05&0.000108456&0.000106767\\
\hline
\end{tabular}
\end{center}
\caption{\label{rate_table} $L^2$ norm of the error for $\eta,u,v$.}
\end{table}
\begin{figure}[h!]
\begin{center}
\includegraphics[height=6cm]{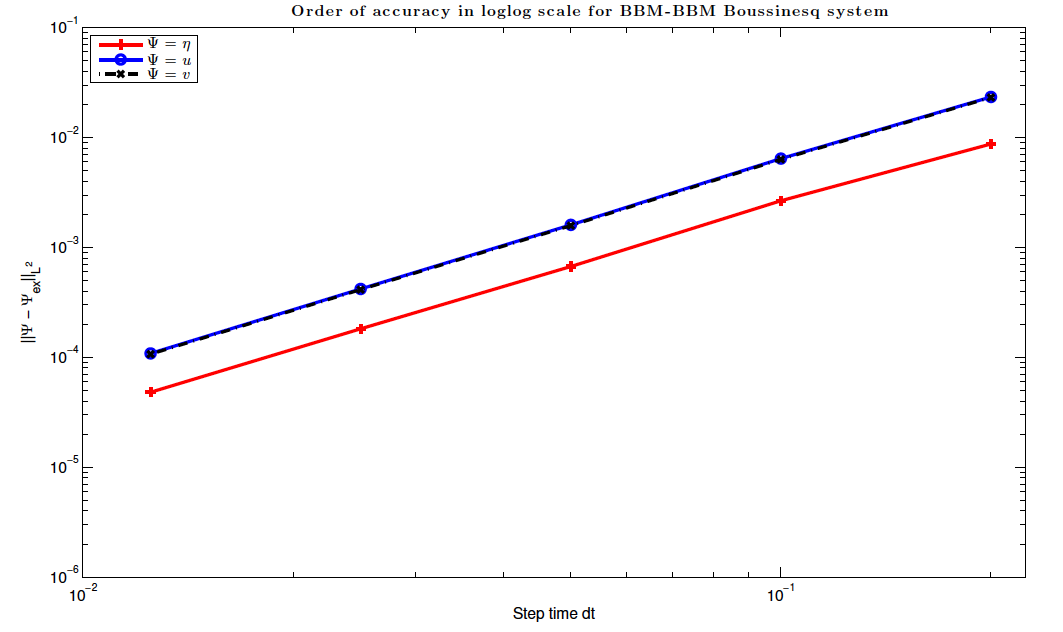}
\caption{\label{rate_BBM} Rate of convergence for BBM-BBM.}
\end{center}
\end{figure}
\subsubsection{Computation time}\label{secopt}
In this section, we consider the BBM-BBM Boussinesq system on the domain $\Omega= ]-40,40[\times ]-40,40[$ with $\mathbb{P}_1$ continuous piecewise linear functions,  the space discretization $\Delta x=0.5$, the time step $\Delta t=0.1$ and as initial data $\eta_{h0}(x,y) = 0.2e^{-(x^2 + y^2)/5},  u_{h0}(x,y) =  v_{h0}(x,y) = 0$ with zero Dirichlet homogenous Boundary Conditions for $\eta_{h}$, $ u_{h}$ and $ v_{h}$ on the whole boundary.\\
In order to solve this system, we will show the time comparison of different method:
\begin{itemize}
\item \textbf{M1} to solve \textbf{Algorithm 1} without using adaptmesh technique and the keyword \textcolor{red}{\ttfamily{init}}.
\item\textbf{M1init} to solve \textbf{Algorithm 1} using the keyword \textcolor{red}{\ttfamily{init}} and without using adaptmesh technique.
\item\textbf{M1A-4} to solve \textbf{Algorithm 1} using adaptmesh technique with {\ttfamily{err=1e-4}} and without the keyword \textcolor{red}{\ttfamily{init}}.
\item\textbf{M1A-2} to solve \textbf{Algorithm 1} using adaptmesh technique with {\ttfamily{err=1e-2}} and without the keyword \textcolor{red}{\ttfamily{init}}.
\item \textbf{M2} to solve \textbf{Algorithm 2} without using adaptmesh technique and the keyword \textcolor{red}{\ttfamily{init}}.
\item\textbf{M2init} to solve \textbf{Algorithm 2} using the keyword \textcolor{red}{\ttfamily{init}} and without using adaptmesh technique.
\item\textbf{M2A-4} to solve \textbf{Algorithm 2} using adaptmesh technique with {\ttfamily{err=1e-4}} and without the keyword \textcolor{red}{\ttfamily{init}}.
\item\textbf{M2A-2} to solve \textbf{Algorithm 2} using adaptmesh technique with {\ttfamily{err=1e-2}} and without the keyword \textcolor{red}{\ttfamily{init}}.
\end{itemize}
We present in Table \ref{adapt_compute}, the time of computation in second at time $T=10s$ using all the different method cited before. All computation was made on a Macbook OS X, Intel core 2 Duo (CPU), 4Go (Memory), 2 Ghz (Processor). 

\begin{center}
\begin{figure}[!htb]
        \includegraphics[width=12cm]{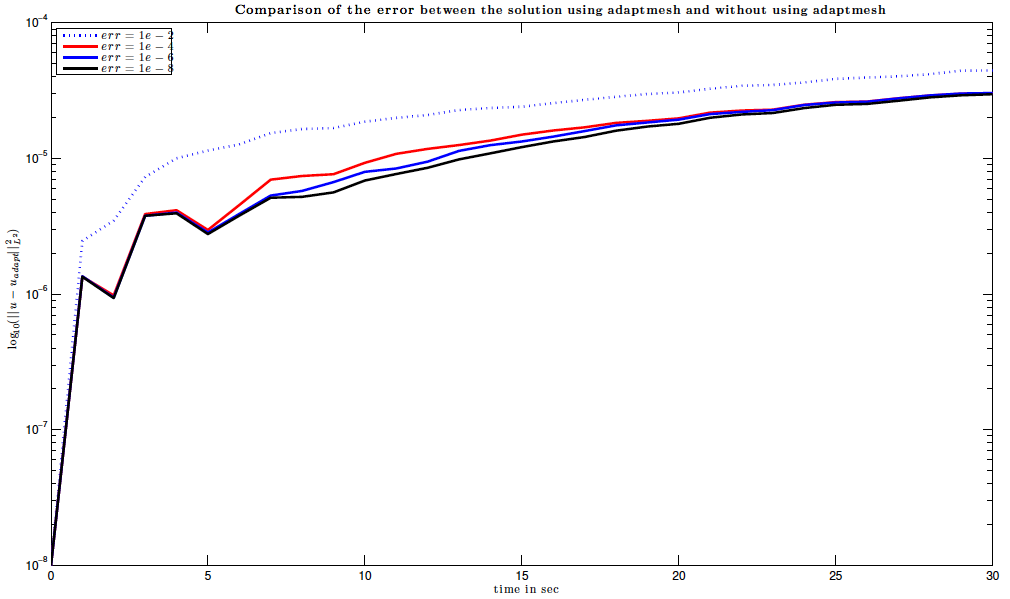}
\caption{\label{comp_adapt}Comparison of the error between the solutions.}
\end{figure}
\end{center}

\begin{center}
\begin{figure}[!htb]
        \includegraphics[width=12cm]{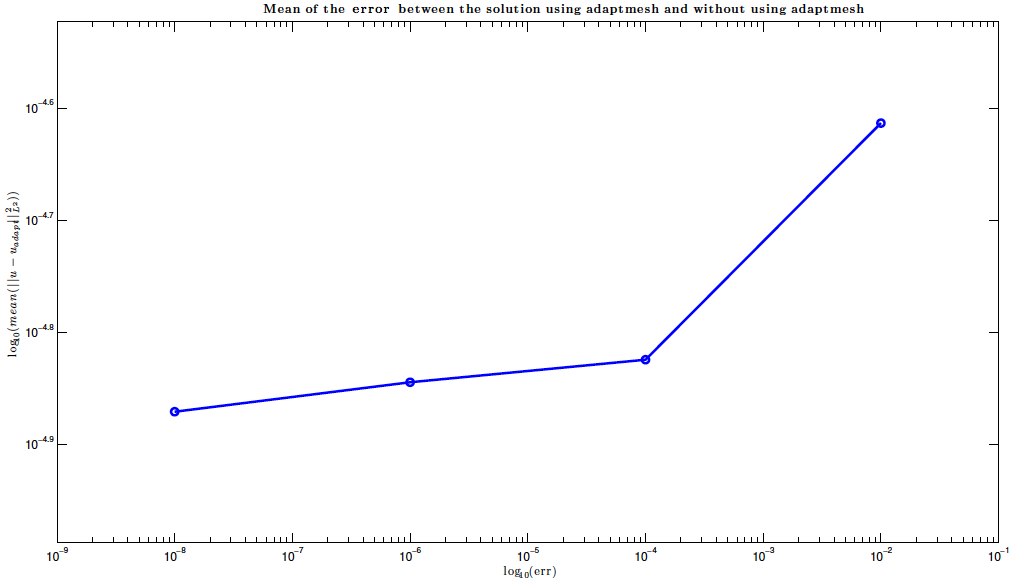}
\caption{\label{mean_comp} Mean of the error between the solutions.}
\end{figure}
\end{center}

\begin{table}[htdp]
\begin{center}
\begin{tabular}{|c|c|c|c|}
\hline
\textbf{M1}&\textbf{M1init}&\textbf{M1A-4}&\textbf{M1A-2}\tabularnewline\hline
1555.48&722.507&115.485&45.3462\tabularnewline\hline
\textbf{M2}&\textbf{M2init}&\textbf{M2A-4}&\textbf{M2A-2}\tabularnewline\hline
1217.01&619.462&99.5689&39.3365\tabularnewline\hline
\end{tabular}
\end{center}
\caption{\label{adapt_compute} Comparison of computation time for the different method used to solve the BBM-BBM system.}
\end{table}
\begin{flushleft}
We note that, without using the mesh adaptation technique, we have the same result for all the computed solution, so we can see from Table \ref{adapt_compute} that the best method to use is the \textbf{M2init}.\vspace{.5cm}\\
In other hand, using the mesh adaptation technique, we can remark from Table \ref{adapt_compute} that the computation time is better then other method, but unfortunately, we have a little difference between the computed solution, that we plot in Figure \ref{comp_adapt} the square of $L^2$ norm of the error between the computed solution using the \textbf{M1} method and the one computed with the \textbf{M1A-8}, \textbf{M1A-6}, \textbf{M1A-4}, \textbf{M1A-2} methods where we have {\ttfamily{err=1e-8}}, {\ttfamily{err=1e-6}}, {\ttfamily{err=1e-4}}, {\ttfamily{err=1e-2}}, respectively vs the time till $T=30s$. We also plot in Figure \ref{mean_comp} the mean of all the square of $L^2$ norm of the error computed for different mesh adaptation method.
\end{flushleft}
\vspace{-1.5cm}
We can remark from this result that we have the same result using \textbf{M1A-8}, \textbf{M1A-6}, \textbf{M1A-4} method and the mean error between the solution computed with these method and the computed one using the \textbf{M1} is of order $10^{-5}$ and we can see the large time difference.\\
\subsubsection{Reflection of expanding symmetric waves at two boundaries of the BBM-BBM Boussinesq system}
\begin{figure}[!htb]
\begin{tabular}{>{\centering}m{6cm}>{\centering}m{6cm}}
        \includegraphics[height=4cm,width=4cm]{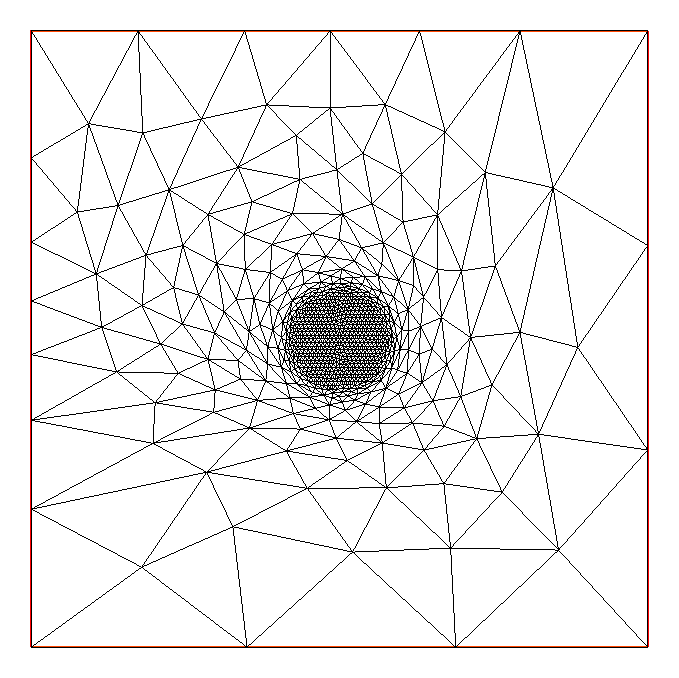}&
        \includegraphics[height=4cm,width=4cm]{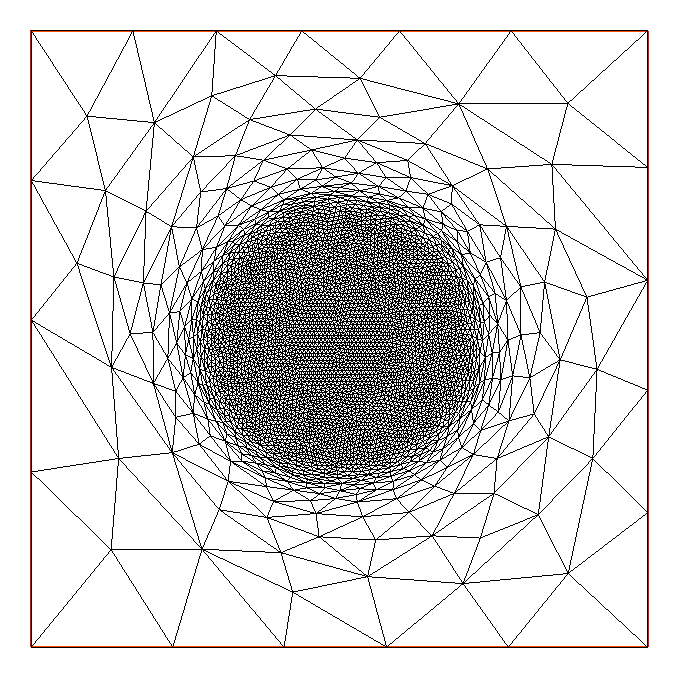}\tabularnewline
        \includegraphics[height=4cm,width=4cm]{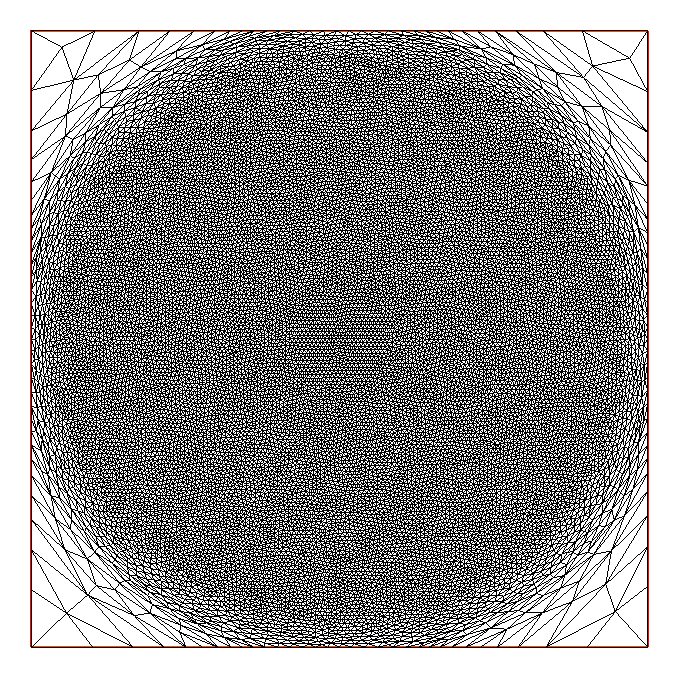}&
        \includegraphics[height=4cm,width=4cm]{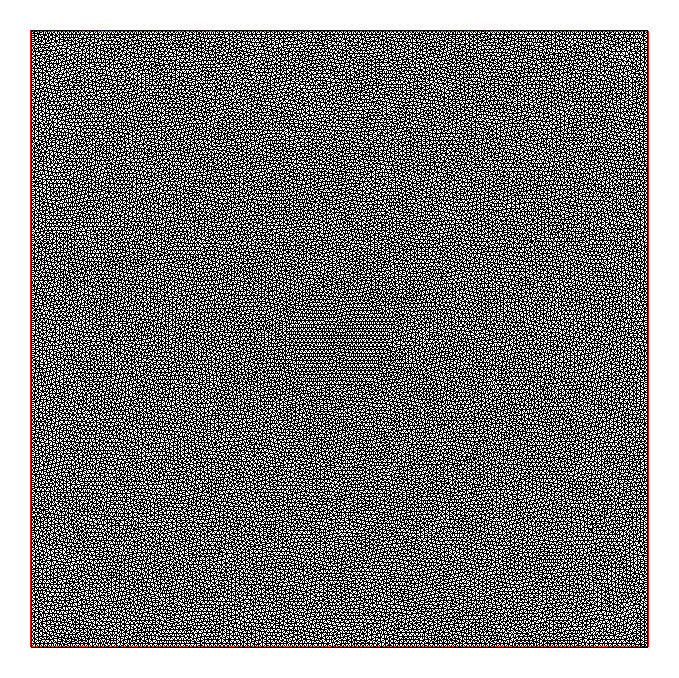}
\end{tabular}
\caption{\label{adaptmesh} Propagation of the mesh for the BBM-BBM Boussinesq system where $a=c=0$ and $b=d=1/6$ for different time $t=\{0.1, 20, 40, 70\}$}
\end{figure}

\begin{figure}[!htb]
\begin{tabular}{>{\centering}m{6cm}>{\centering}m{6cm}}
        \includegraphics[height=4cm,width=5cm]{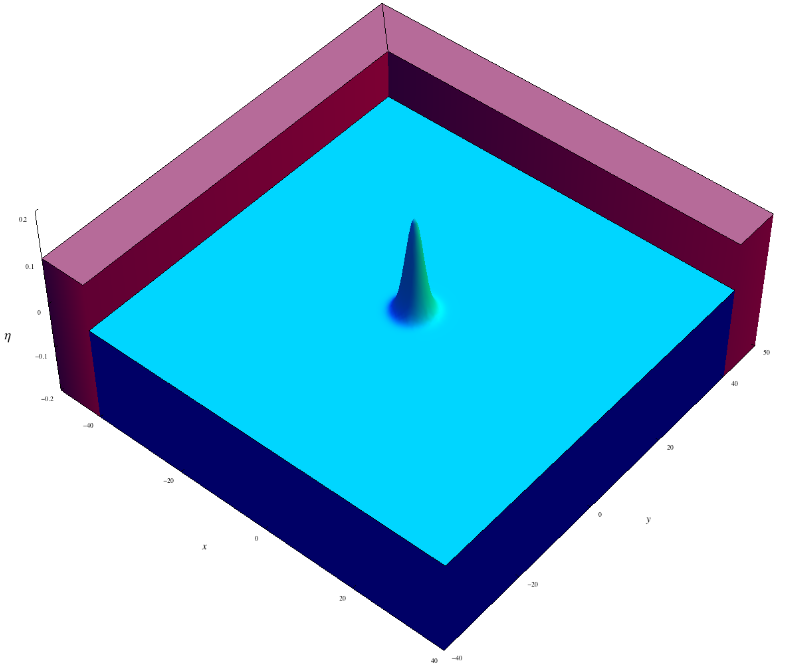}&
        \includegraphics[height=4cm,width=5cm]{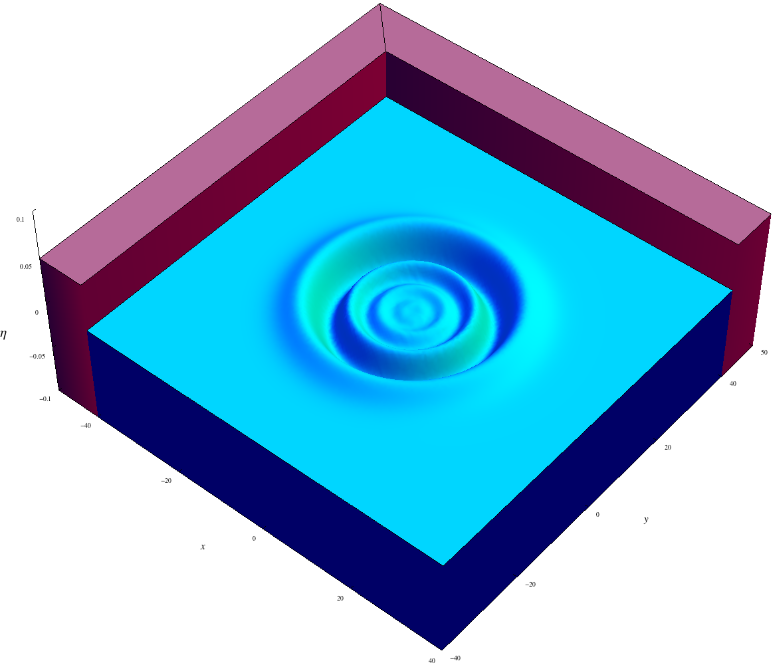} \tabularnewline
        \includegraphics[height=4cm,width=5cm]{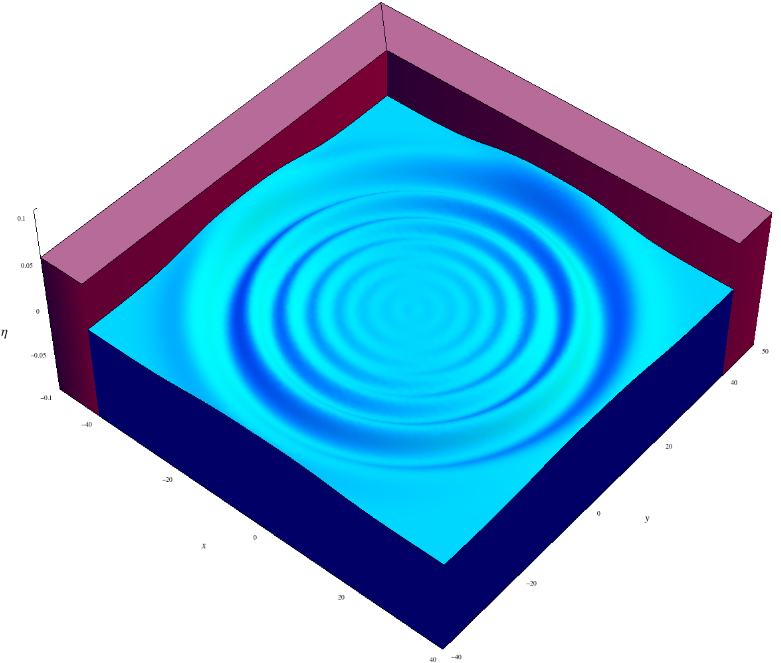}&
        \includegraphics[height=4cm,width=5cm]{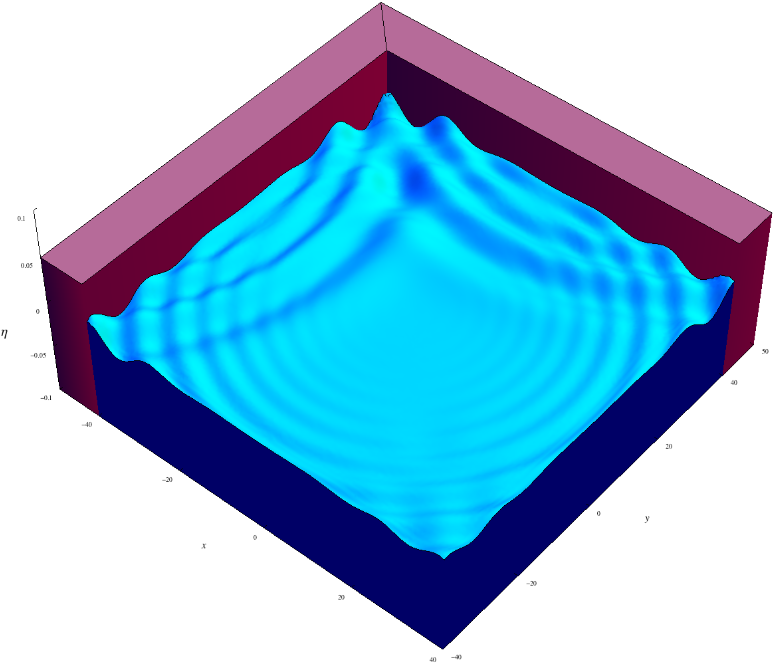} 
\end{tabular}
\caption{\label{BBM} Propagation of the solution of the BBM-BBM Boussinesq system where $a=c=0$ and $b=d=1/6$ for different time $t=\{0, 20, 40, 70\}$}
\end{figure}

In Figures \ref{adaptmesh} and \ref{BBM} we show the reflection from two parts of the boundary of an expanding symmetric wave of the BBM-BBM Boussinesq system where  $a=c=0$ and $b=d=1/6$. For this experiment we used as initial data the functions $\eta_{h0}(x,y) = .2e^{-(x^2 + y^2)/5},  u_{h0}(x,y) =  v_{h0}(x,y) = 0$. We used zero Neumann Boundary Conditions for $\eta_{h}$ on the whole boundary, zero Dirichlet data for $ u_{h}$ and $ v_{h}$ on  $x =-40$ and $y = 40$ (where we have the wall), and zero Neumann boundary data for $ u_{h}$ and $ v_{h}$ on $x = 40$ and $y = -40$. The expanding waves are reflected from the $x = -40$ and $y = 40$ parts of the boundary.\newline
We note that in Figure \ref{adaptmesh} we show the effect of the mesh adaptation following the evolution of $\eta_{h}$ in time and in Figure \ref{BBM} we show the propagation of the solution $\eta_h$.
\clearpage
\subsubsection{Expanding symmetric waves under the KdV-KdV Boussinesq system}
\begin{figure}[!h]
\begin{tabular}{>{\centering}m{6cm}>{\centering}m{6cm}}
        \includegraphics[height=4cm,width=5cm]{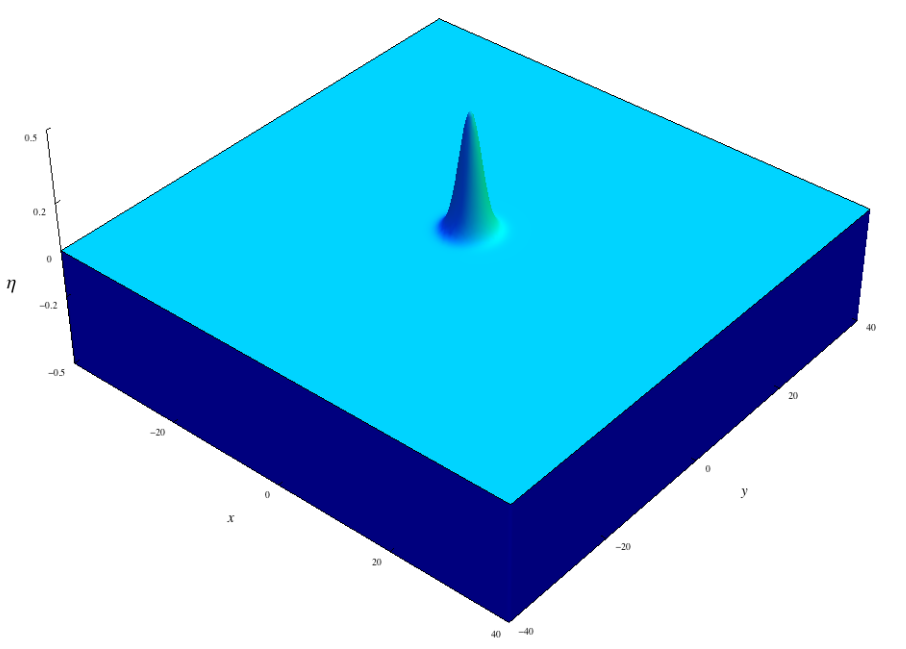}&
        \includegraphics[height=4cm,width=5cm]{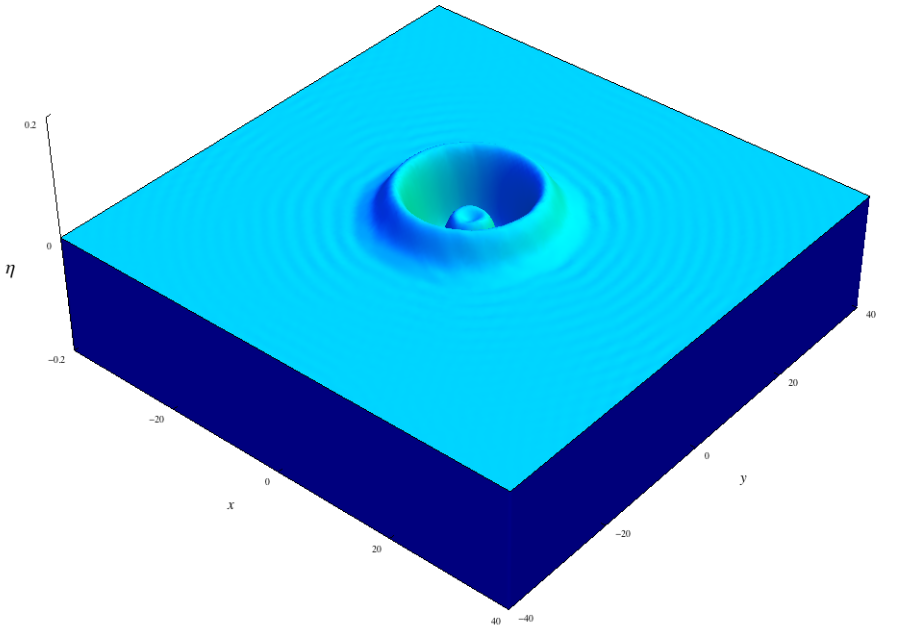} \tabularnewline
        \includegraphics[height=4cm,width=5cm]{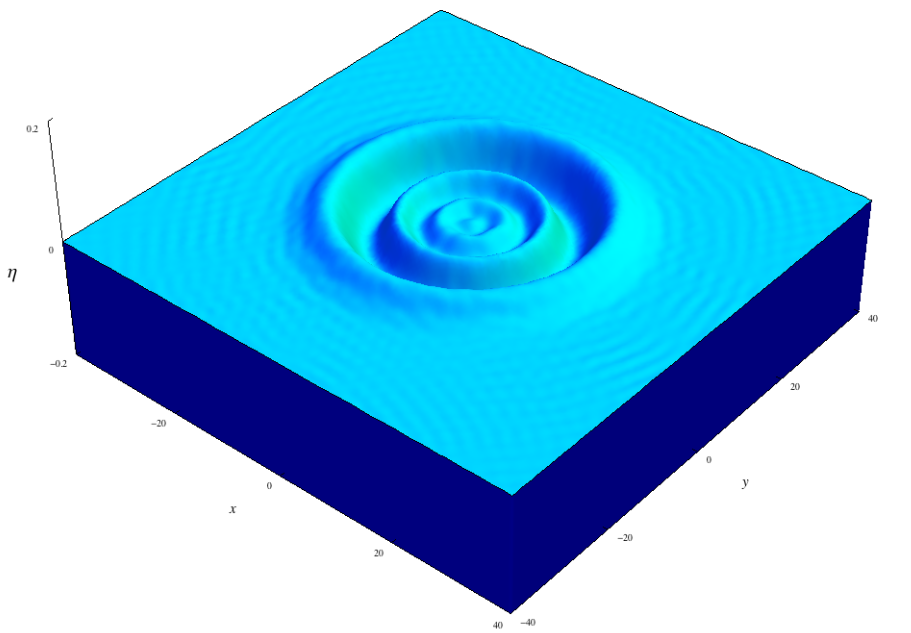} &
        \includegraphics[height=4cm,width=5cm]{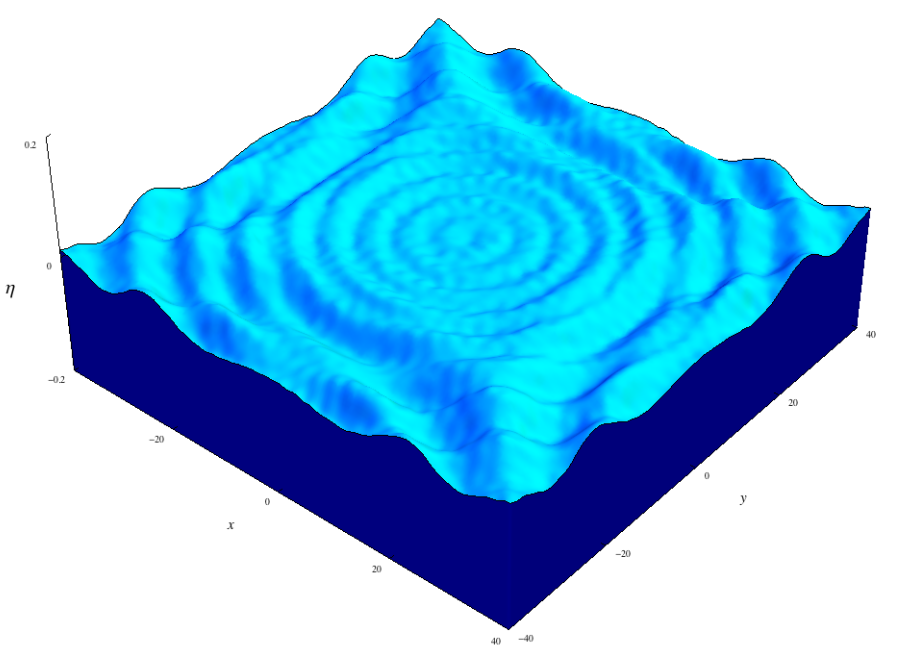} 
\end{tabular}
\caption{\label{KdV-KdV} Propagation of the solution of the KdV-KdV Boussinesq system where $a=c=1/6$ and $b=d=0$ for different time $t=\{0,10,20,60\}$}
\end{figure}

In Figure \ref{KdV-KdV}, we present the evolution of the $\eta_h$ profile emanating from the radially symmetric initial data $\eta_{h0}(x,y) = .5e^{-(x^2 + y^2)/5},  u_{h0}(x,y) =  v_{h0}(x,y) = 0$, under the KdV-KdV Boussinesq system where  $a=c=1/6$ and $b=d=0$. We used Bi-Periodic Boundary Conditions for $\eta_{h}$, $ u_{h}$ and $ v_{h}$ and we work with the time step $\Delta t=0.001$.
We remark here that with these Bi-Periodic Boundary Conditions for $\eta$, $ u$ and $ v$ and their derivatives, in addition by integrating the equations in the system (\ref{BOUS}) on the hole domaine, we deduce the following mass conservation:
$(Id-b\Delta)\ds\int_{\Omega}\eta_t=0$ and the relations $(Id-d\Delta)\ds\int_{\Omega} u_t=0$, $(Id-d\Delta)\ds\int_{\Omega} v_t=0$. Hence: 
\bg\label{CONSMASS}
\int_{\Omega}\eta=cte=\int_{\Omega}\eta_0, \qquad \int_{\Omega} u=cte=\int_{\Omega} u_0, \qquad \int_{\Omega} v=cte=\int_{\Omega} v_0.\ed
In other hand, numerically, we see that these defined quantity are well conserved over time and we have: 
$$\int_{\Omega}\eta=cte=\int_{\Omega}\eta_0=7.84527, \qquad \int_{\Omega} u=cte=\int_{\Omega} u_0=0, \qquad \int_{\Omega} v=cte=\int_{\Omega} v_0=0.$$
\vspace{-.5cm}
\begin{flushleft}
We can see in Figure \ref{KdV-KdV} from $t=10$ a small amplitude periodic profile (ripples) which are propagating in front of the wavefront and which has been observed in \cite{Bona*07}\label{Bona*071} (in the case of 1D KdV-KdV system). These ripples are still observed when using $\mathbb{P}_2$ elements but they do not infer to the numerical stability, we still have the mass conservation of the fluid.\\
Other simulation for different Boussinesq systems can be found in \cite{Sad11}\label{Sad112}.
\end{flushleft}
\subsection{Comparison of the results with Mitsotakis {\it{et al.}}}
\begin{figure}[!htb]
\begin{tabular}{>{\centering}m{6cm}>{\centering}m{6cm}}
\includegraphics[height=6cm,width=6cm]{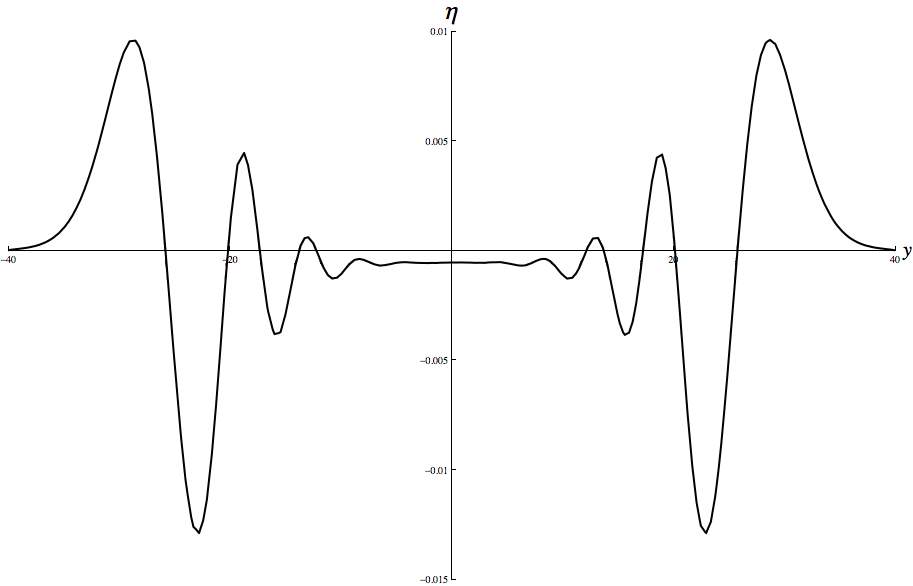} &
\includegraphics[height=6cm,width=6cm]{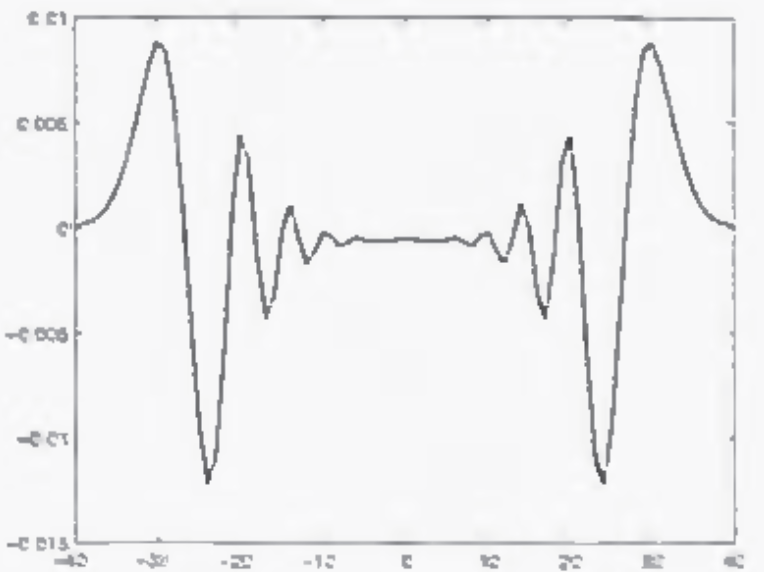} \tabularnewline
\includegraphics[height=6cm,width=6cm]{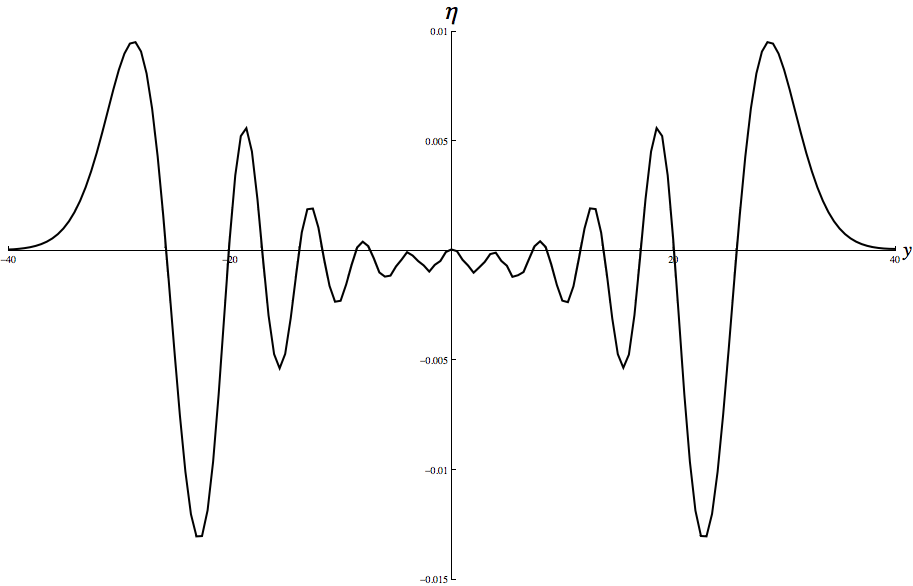}&
\includegraphics[height=6cm,width=6cm]{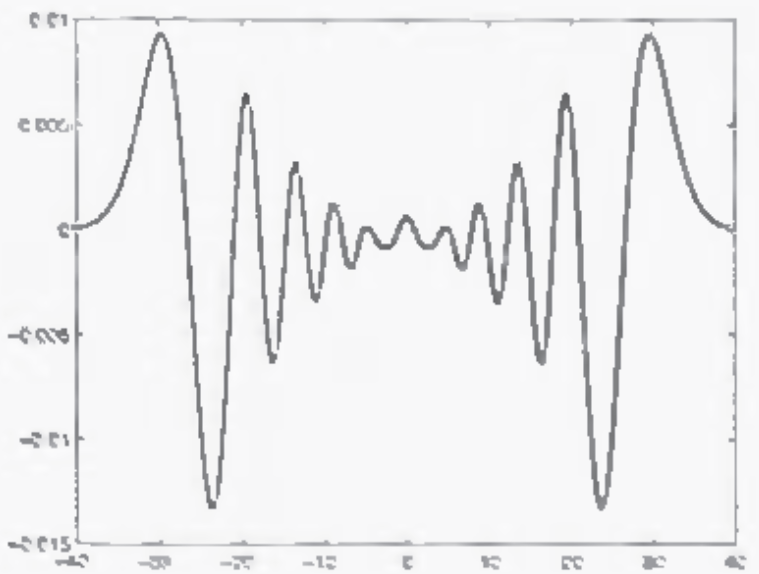} 
\end{tabular}
\caption{\label{COMPDM} Cross sections in $x$ direction of $\eta(x, y,t)$ at $t = 30s$, Bona-Smith (above) and BBM-BBM (below), using our code (left) and Mitsotakis {\it et al.} code (right) borrowed from the article \cite{Dou*07}\label{Dou*079}.}
\end{figure}

We validate also our code by comparing the results with those of Mitsotakis {\it{et al.}} \cite{Dou*07} (in Figure 8, page 847) where they consider the cross section in $x$ direction to the $\eta$ component of the solution from radially symmetric initial data of the form $\eta_{h0}(x,y) = .2e^{-(x^2 + y^2)/5}$ with $u_{h0}(x,y) =  v_{h0}(x,y) = 0$ under the BBM-BBM and the Bona-Smith system, both considered with zero Dirichlet boundary conditions for $\eta,u$ and $v$. We remark in Figure \ref{COMPDM}, that the shape of the solution for the results of our code at the left part and for those of Mitsotakis {\it et al.} code at the right part are similar.\\
We note also that in their paper, they used Galerkin finite elements based on tensor products of smooth splines and an explicit second order (for BBM-BBM system with bilinear splines) and fourth order (for Bona-Smith system with bi-cubic splines) Runge-Kutta scheme while in the present work, we used $\mathbb{P}_1$ continuous piecewise linear functions and an explicit second order Runge-Kutta scheme for all systems.\\

\subsection{Comparison of different Boussinesq models}
\begin{figure}[!htb]
\begin{tabular}{>{\centering}m{6cm}>{\centering}m{6cm}}
\includegraphics[height=4cm,width=6cm]{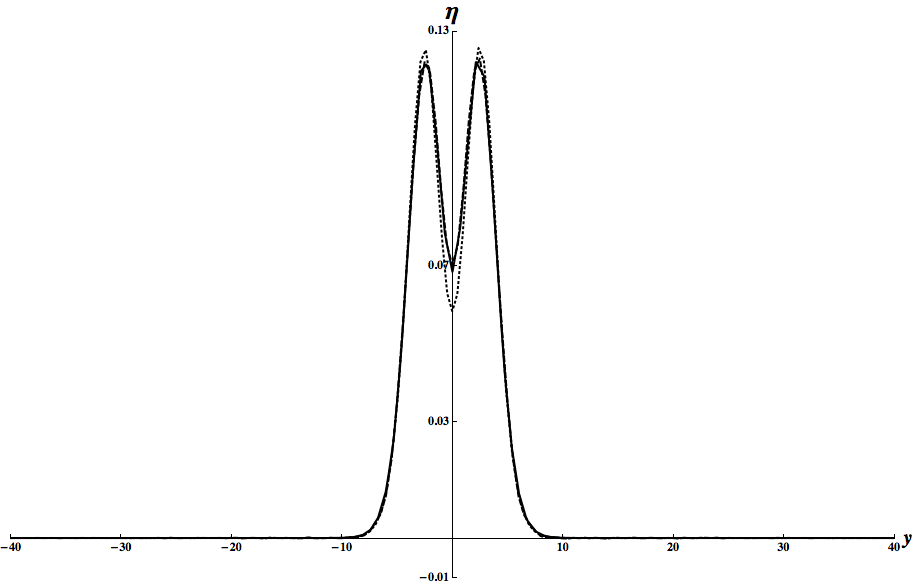} &
\includegraphics[height=4cm,width=6cm]{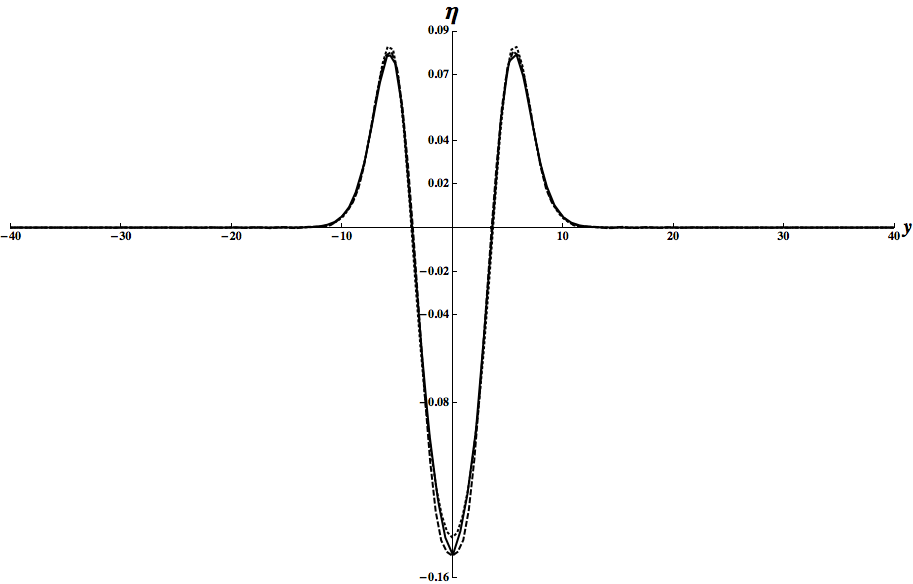} \tabularnewline
\includegraphics[height=4cm,width=6cm]{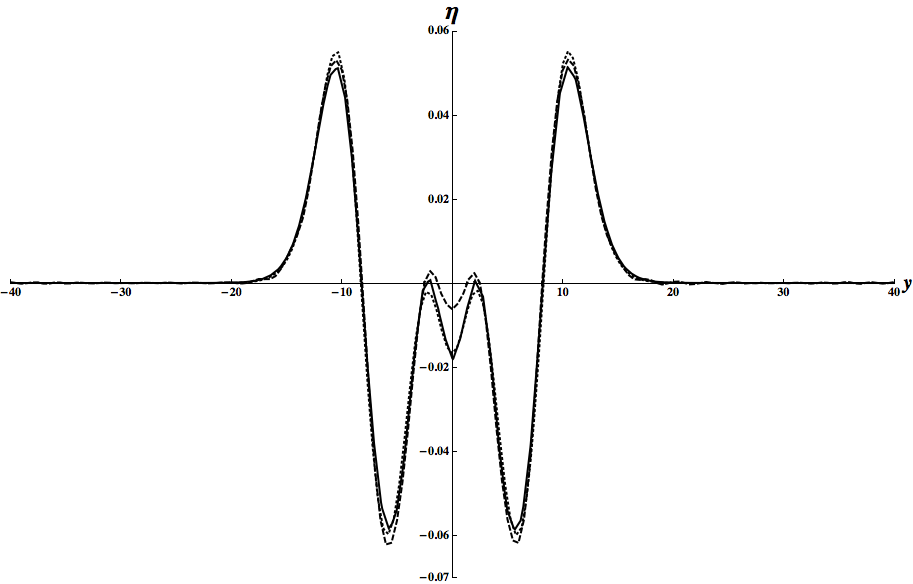}&
\includegraphics[height=4cm,width=6cm]{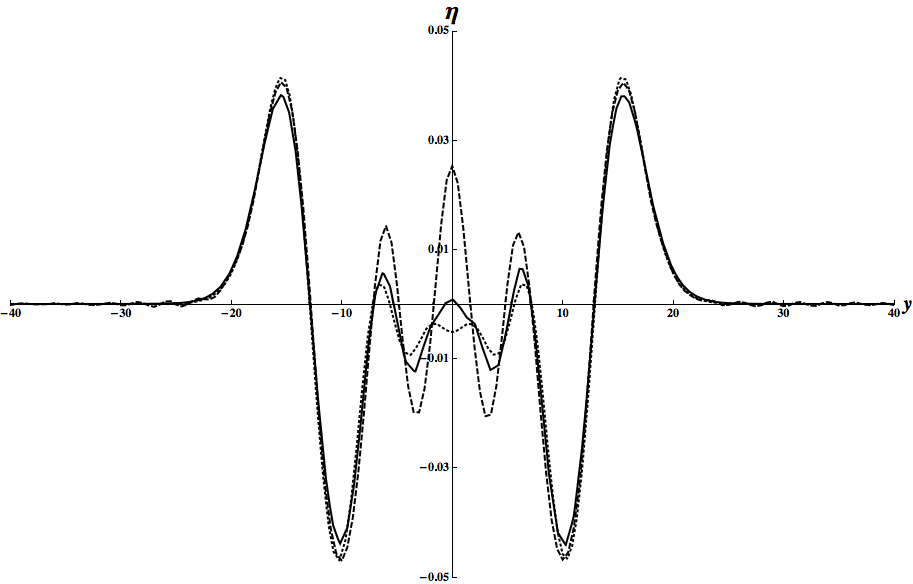} \tabularnewline
\includegraphics[height=4cm,width=6cm]{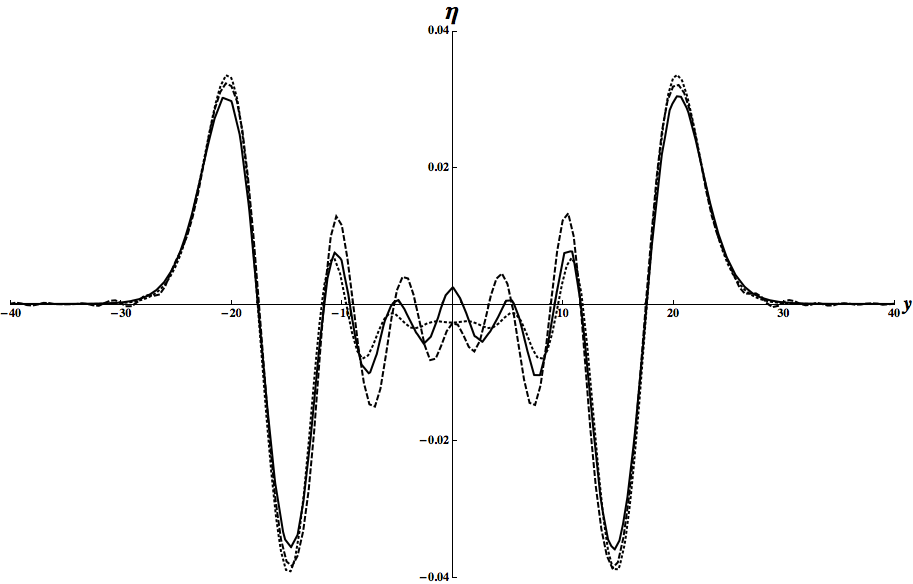}&
\includegraphics[height=4cm,width=6cm]{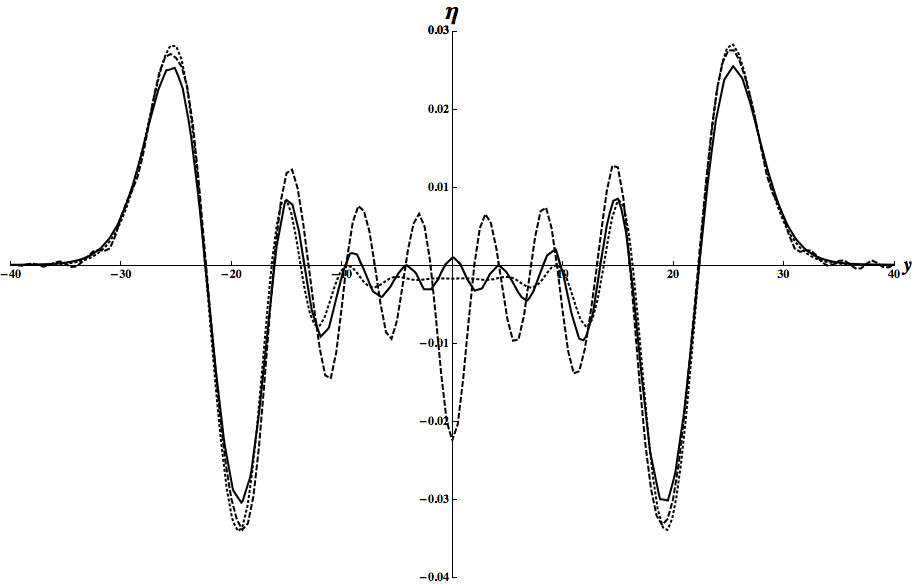} 
\end{tabular}
\caption{\label{COMPKBS} Comparison of the cross sections in $x$ direction of $\eta(x, y,t)$ at time $t = \{2,5,10,15,20,25\}s$, where in all the figure, the dotted line is for Bona-Smith system, the dashed for KdV-KdV one and the solid line for BBM-BBM one.}
\end{figure}

We compare here KdV-KdV, BBM-BBM and Bona-Smith models as defined in section \ref{probset} Table \ref{TABBOUS2D}.\\
In Figures \ref{COMPKBS}, we present a comparison of the evolution of the $\eta$ component of the solution from radially symmetric initial data of the form $\eta_{h0}(x,y) = .5e^{-(x^2 + y^2)/5}$ with $u_{h0}(x,y) =  v_{h0}(x,y) = 0$, under the BBM-BBM
(solid line with zero Dirichlet b.c., $\Delta t=0.1$, $\Delta x=0.5$), the Bona-Smith system (dotted line with zero Dirichlet b.c., $\Delta t=0.1$, $\Delta x=0.5$) and the KdV-KdV system (dashed line with periodic b.c., $\Delta t=0.001$, $\Delta x=0.5$). Figure \ref{COMPKBS} shows the cross sections of the $\eta$ profiles for different time in the $x$ - direction. The speed and the amplitude of the outgoing front is approximately the same for the BBM-BBM and Bona-Smith systems but the pattern of the oscillations behind the fronts are different: in the case of the Bona-Smith system the two outgoing wave trains have practically separated by $t = 25s$, while the larger in amplitude dispersive oscillatory tails of the BBM-BBM solution seem to be still interacting. In other hand, the speed of the outgoing front for the KdV-KdV system is approximately the same with other systems while the amplitude for the internal crest is bigger, we remark also that the solution still interacting after $t=25s$.

\section{Conclusion}
We have presented a numerical approach with \freefem to solve the Boussinesq systems with a flat bottom, we validated our code and establishes the adequacy of the chosen finite element discretization by comparing the results with those of Mitsotakis {\it et al.} We have established also the feasibility of simulating complex equations of hydrodynamics as Boussinesq systems with \freefem and we have optimized the algorithm that we use in section \ref{secopt}.\\
Using this approach, we can consider the case of a variable bottom (in space and/or in time), see \cite{Sad11}\label{Sad113} which is an ongoing work to appear soon. As a feature we address the simulation of {\it Tsunamis} with our approach, by including realistic data (bathymetry, generation of tsunami waves).\vspace{1cm}\\
{\bf Acknowledgements: }This work is supported by the regional program "Appui \`a lÕ\'emergence" of the Region Picardie. I would like to thank my PhD advisor Jean-Paul Chehab (LAMFA, Amiens), Denys Dutykh (LAMA, Savoie), Gloria Faccanoni (IMATH, Toulon), Kirill Pichon Gostaf (LJLL, Paris), Fr\'ed\'eric Hecht (LJLL, Paris), Antoine Le Hyaric (LJLL, Paris), Youcef Mammeri (LAMFA, Amiens), Olivier Pantz (CMAP, Paris) and Dimitrios Mitsotakis (IMA, University of Minnesota) for fruitful discussions and remarks.
\clearpage
\bibliographystyle{plain}

\begin{thebibliography}{}

\end{thebibliography}


\begin{thebibliography}{99}
\bibitem{Aze*080}{\sc Pascal Azerad, Frederic Bouchette, Benjamin Ivorra, Damien Isebe and 
 Mohammadi}. {\sl Shape optimization of geotextile tubes for sandy beach protection. } \href{http://www.math.univ-montp2.fr/~azerad/hdr/ijnmeonline.pdf}{{\em Int. J. Numer. Meth. Engng}, Vol. 74, pp. 1262-1277, 2008. }\\

\bibitem{Aze**080}{\sc Pascal Azerad, Frederic Bouchette, Damien Isebe and Bijan Mohammadi}. {\sl Optimal shape design of defense structures for minimizing short wave impact. } \href{http://www.sciencedirect.com/science?_ob=MiamiImageURL&_cid=271917&_user=781134&_pii=S0378383907000725&_check=y&_origin=&_coverDate=31-Jan-2008&view=c&wchp=dGLzVlB-zSkWA&md5=e530525704be89082068fad23366effb/1-s2.0-S0378383907000725-main.pdf}{{\em Coastal Engineering}, Vol. 55, pp. 35-46, 2008. }\\

\bibitem{Bona*04} {\sc Jerry Lioyd Bona, Min Chen and Jean-Claude Saut}. {\sl Boussinesq equations and other systems for small amplitude long waves in nonlinear dispersive media: II. the nonlinear theory.} \href{http://www.math.uic.edu/~bona/papers/boussineq-paper.pdf}{{\em Nonlinearity}, 17, no3, 925-952, 2004. }\\

\bibitem{Bona*07} {\sc J. L. Bona, V. A. Dougalis and D. E. Mitsotakis}. {\sl Numerical solution of KdV-KdV systems of Boussinesq equations I: The numerical scheme and generalized solitary waves}. \href{http://sites.google.com/site/dmitsot/BDM1.pdf?attredirects=0}{{\em Math. Comput. Simulation}, 74:214-228, 2007.}\\

\bibitem{Bre*94}{\sc Suzanne C. Brenner and L. Ridgway Scott}. {\sl The Mathematical Theory of Finite Element Methods.} \href{http://www.springerlink.com/content/978-0-387-75934-0\#section=144131&page=1}{{\em Springer-Verlag}, New York, 1994. }\\

\bibitem{Chen09} {\sc Min Chen}. {\sl Numerical investigation of a two-dimensional Boussinesq system.} \href{http://www.math.purdue.edu/~chen/papers/tsunami.pdf}{{\em Discrete and Continuous Dynamical Systems}, Vol 23, no4, 1169-1190, April 2009. }\\

\bibitem{Chen*09} {\sc Min Chen and Olivier Goubet}. {\sl Long-time asymptotic behavior of two-dimensional dissipative Boussinesq systems. }\href{http://www.math.purdue.edu/~chen/papers/goubet2.pdf}{{\em Discrete and Continuous Dynamical Systems Series S}, Vol 2, no1, 37-53, March 2009. }\\

\bibitem{Dem91} {\sc Jean-Pierre Demailly}. {\sl Analyse num\'erique et \'equations diff\'erentielles. } \href{http://grenoble-sciences.ujf-grenoble.fr/ouvrages/extraits/13011171.pdf}{{\em Presses Universitaires de grenoble}, 1991. }\\

\bibitem{Dou*07} {\sc Vassilios Dougalis, Dimitrios Mitsotakis and Jean-Claude Saut}. {\sl On some Boussinesq systems in two space dimensions: theory and numerical analysis. } \href{http://sites.google.com/site/dmitsot/dms1.pdf?attredirects=0}{{\em M2AN}, no. 5, 825-854, Vol 41, 2007. }\\

\bibitem{Dou*09} {\sc Vassilios Dougalis, Dimitrios Mitsotakis and Jean-Claude Saut}. {\sl On initial-boundary value problems for a Boussinesq system of BBM-BBM type in a plane domain. } \href{http://sites.google.com/site/dmitsot/dms2.pdf?attredirects=0}{ {\em AIMS' Journal}, Vol 23, 2009. }\\

\bibitem{Dou*10} {\sc Vassilios Dougalis, Dimitrios Mitsotakis and Jean-Claude Saut}. {\sl Boussinesq systems of Bona-Smith type on plane domain: Theory and Numerical Analysis.} \href{http://sites.google.com/site/dmitsot/dms3.pdf?attredirects=0}{{\em Journal of Scientific Computing}, Vol. 44, no2, pp. 109-135, 2010. }\\

\bibitem{Dut*07} {\sc Denys Dutykh and Fr\'ed\'eric Dias}. {\sl Dissipative Boussinesq equations.} \href{http://hal.archives-ouvertes.fr/docs/00/18/58/20/PDF/dissbouss.pdf}{ {\em C. R. M\'ecanique}, 335 , 559-583, 2007. }\\

\bibitem{Dut*110} {\sc Denys Dutykh, Theodoros Katsaounis and Dimitrios Mitsotakis.} {\sl Finite volume schemes for dispersive wave propagation and runup. } \href{http://hal.archives-ouvertes.fr/docs/00/55/37/43/PDF/DKM-2010.pdf}{ {\em Computational Physics}, 230 (8), 3035 - 3061, 2011. }\\

\bibitem{Lin*113}{\sc F\'elipe Linares, Didier Pilod and Jean-Claude Saut.} {\sl Well-posedness of strongly dispersive two-dimensional surface waves Boussinesq.} \href{http://www.im.ufrj.br/~didier/Research/LPS(15).pdf}{{\em ArXiv}:1103.4159v2, 11 Apr 2011. }\\

\bibitem{Luc*98} {\sc Brigitte Lucquin and Olivier Pironneau}. {\sl Introduction to Scientific Computing}. \href{http://www.abebooks.com/9780471972662/Introduction-Scientific-Computing-Lucquin-Brigitte-0471972665/plp}{ {\em Wiley}, 1998. }\\

\bibitem{Mir*05} {\sc Alain Miranville and Roger Temam}. {\sl Mathematical modeling in continuum mechanics.} \href{http://ebooks.cambridge.org/ebook.jsf?bid=CBO9780511755422}{{\em Cambridge University Press}, 2005.}\\

\bibitem{Mit09} {\sc Dimitrios Mitsotakis}. {\sl Boussinesq systems in two space dimensions over a variable bottom for the generation and propagation of tsunami waves.} \href{http://sites.google.com/site/dmitsot/tsunami_2007.pdf?attredirects=0}{{\em Mat. Comp. Simul.}, 80:860-873, 2009.}\\

\bibitem{Sad11} {\sc Georges Sadaka}. {\sl Etude math\'ematique et num\'erique d'\'equations d'ondes aquatiques amorties}.  \href{http://lamfa.u-picardie.fr/sadaka/these_GS.pdf}{{\em Th\`ese de l'Universit\'e de Picardie Jules Verne - Amiens}, 2011.}\\

\bibitem{Wal*02} {\sc Mark Walkley and Martin Berzins}. {\sl A finite element method for the two-dimensional extended Boussinesq equations. } \href{http://onlinelibrary.wiley.com/doi/10.1002/fld.349/pdf}{{\em Int. J. Numer. Meth. Fluids}, 39:865-885, 2002}.\\

\end{thebibliography}

\end{document}